# A POSTERIORI ERROR ESTIMATES FOR THE FINITE ELEMENT DISCRETIZATION OF SECOND-ORDER PDES SET IN UNBOUNDED DOMAINS

T. CHAUMONT-FRELET*

Abstract. We consider second-order PDE problems set in unbounded domains and discretized by Lagrange finite elements on a finite mesh, thus introducing an artificial boundary in the discretization. Specifically, we consider the reaction diffusion equation as well as Helmholtz problems in waveguides with perfectly matched layers. The usual procedure to deal with such problems is to first consider a modeling error due to the introduction of the artificial boundary, and estimate the remaining discretization error with a standard a posteriori technique. A shortcoming of this method, however, is that it is typically hard to obtain sharp bounds on the modeling error. In this work, we propose a new technique that allows to control the whole error by an a posteriori error estimator. Specifically, we propose a flux-equilibrated estimator that is slightly modified to handle the truncation boundary. For the reaction diffusion equation, we obtain fully-computable guaranteed error bounds, and the estimator is locally efficient and polynomial-degree-robust provided that the elements touching the truncation boundary are not too refined. This last condition may be seen as an extension of the notion of shape-regularity of the mesh, and does not prevent the design of efficient adaptive algorithms. For the Helmholtz problem, as usual, these statements remain valid if the mesh is sufficiently refined. Our theoretical findings are completed with numerical examples which indicate that the estimator is suited to drive optimal adaptive mesh refinements.

## 1. Introduction

Partial differential equations (PDE) set in an unbounded domain $\Omega$ occur in many applications, in particular in wave propagation problems. When discretizing such problem with the finite element method, it is common to first introduce an approximate PDE model set on a bounded domain $\widetilde{\Omega} \subset \Omega$, and discretize the resulting problem on a mesh of $\widetilde{\Omega}$, see e.g. [21, 24]. The error analysis then takes the form

$$e := \|u - \widetilde{u}_h\| \leq \|u - \widetilde{u}\| + \|\widetilde{u} - \widetilde{u}_h\| =: \widetilde{\varepsilon} + \varepsilon_h,$$

where $u$ is the solution of the unbounded problem, $\widetilde{u}$ the solution to the PDE problem on $\widetilde{\Omega}$ and $u_h$ its finite element discretization. In this conceptual view, the discretization error $\varepsilon_h$ can be estimated with any standard technique, since it concerns a bounded domain. The drawback, however, is that it is typically hard to estimate to modelling error $\widetilde{\varepsilon}$, since it has to rely on PDE analysis. Besides, in the context of adaptive algorithms, it may be hard to balance these two errors to properly set up the truncation boundary in $\widetilde{\Omega}$ and refine the underlying mesh.

*Inria, Univ. Lille, CNRS, UMR 8524 – Laboratoire Paul Painlevé

In this work, we follow an alternative path, where, given a conforming mesh $\mathcal{T}_h$ covering $\widetilde{\Omega}$, we view the finite element space $V_h$ as a subset of $H^1_0(\Omega)$ rather than $H^1_0(\widetilde{\Omega})$. The notation for Sobolev spaces is rigorously introduced in Section 2.4 below. Although this does not change the numerical scheme and the discrete solution, the goal of this viewpoint is to immediately consider the entire error $e$ as a discretization error, without introducing an auxiliary PDE solution $\widetilde{u}$.

To give more precise statements, let us first consider the reaction diffusion equation set in $\Omega \subset \mathbb{R}^d$, $d = 2$ or 3. Given $f \in L^2(\Omega)$, we look for $u \in H^1(\Omega)$ such that

$$\begin{cases} \kappa^2 u - \Delta u &= f \quad \text{in } \Omega,\\ u &= 0 \quad \text{on } \partial\Omega, \end{cases}$$

where the boundary $\partial\Omega$ is allowed to be empty, bounded or even unbounded, as long as it can be exactly meshed. The reaction coefficient $\kappa : \Omega \to \mathbb{R}$ is uniformly bounded away from 0 and $+\infty$. Notice that we do not assume that $f$ is compactly supported nor that $\kappa$ is constant outside a compact set.

Given a conforming simplicial mesh $\mathcal{T}_h$ covering $\Omega^h \subset \Omega$, we consider a $H^1_0(\Omega^h)$-conforming $hp$ finite element space $V_h$, thus introducing an interface $\Gamma^h := \partial\Omega \setminus \partial\Omega^h$ where the discrete solution $u_h \in V_h$ artificially vanishes. For simplicity, in this introduction, we assume that $\kappa$ and $f$ are piecewise constant onto the mesh $\mathcal{T}_h$.

We propose a flux-equilibrated estimator with a flux $\boldsymbol{\sigma}_h \in \boldsymbol{H}(\mathrm{div}, \Omega^h)$ constructed following the standard procedure in [4, 10, 14]. However, the usual estimator $\|\boldsymbol{\nabla} u_h + \boldsymbol{\sigma}_h\|_{\Omega^h}$, cannot be readily employ because both $u \neq 0$ and $\boldsymbol{\sigma}_h \cdot \boldsymbol{n} \neq 0$ on $\Gamma^h$ in general. Hence, we derive the alternative error estimate

$$\|\kappa(u - u_h)\|^2_\Omega + \|\boldsymbol{\nabla}(u - u_h)\|^2_\Omega \leq \sum_{K \in \mathcal{T}_h} \eta_K^2 + \|\kappa^{-1} f\|^2_{\Omega\setminus\Omega^h}, \tag{1.1}$$

with

$$\eta_K := \|\boldsymbol{\sigma}_h + \boldsymbol{\nabla} u_h\|_K + \mu_K \rho_K^{1/2} \|\boldsymbol{\sigma} \cdot \boldsymbol{n}\|_{\partial K \cap \Gamma^h},$$

where $\rho_K$ is the radius of the largest ball contained in $K$ and $\mu_K$ is a fully computable constant scaling as $(\kappa_K h_K)^{-1}$. We can also show that

$$\eta_K \leq C\left\{\left(1 + \delta_{\partial K \cap \Gamma^h} \frac{p_K^+}{\kappa_K h_K}\right)\left(1 + \frac{\kappa_K^+ h_K}{p_K^-}\right)\left(\|\kappa(u - u_h)\|^2_{\widetilde{K}} + \|\boldsymbol{\nabla}(u - u_h)\|^2_{\widetilde{K}}\right)^{1/2} + \delta_{\partial K \cap \Gamma^h} \frac{p_K^+}{\kappa_K h_K} \|\boldsymbol{\nabla} u\|_K\right\}, \tag{1.2}$$

where $\widetilde{K}$ is the element patch around $K$, $\kappa_K^+$ is the maximal value of $\kappa$ in the patch, and similarly $p_K^+$ and $p_K^-$ are the maximal and minimal value of the polynomial degree in the neighboring elements. The symbol $\delta_{\partial K \cap \Gamma^h} = 0$ if $\partial K \cap \Gamma^h$ has zero surface measure and 1 otherwise. The constant $C$ only depends on the shape-regularity parameter of the mesh around $K$, and is independent of polynomial degrees.

The reliability estimate in (1.1) is fully satisfactory. The efficiency estimate in (1.2) is completely standard for elements no touching the artificial boundary, see e.g. [27, 28]. For elements $K$ touching the artificial boundary, however, there are two apparently problematic terms in the estimate, namely $p_K^+/(\kappa_K h_K)$ and $\|\boldsymbol{\nabla} u\|_K$. We claim that neither is an issue.

First, concerning the $p_K^+/(\kappa_K h_K)$ term, it probably makes little sense anyway to use very refined elements or high polynomial degrees touching (and in fact, close to) the artificial boundary. Indeed, the approximation of $u$ by 0 on $\Gamma^h$ is so crude that it seems unlikely to increase the accuracy by local refinements. Actually, optimally refined $hp$ spaces are such that $p_K^+/(\kappa_K h_K)$ is small close to the truncation boundary. In addition, one can easily imagine how to design adaptive mesh refinement algorithms that preserve this property, and we propose one such algorithm in Section 6 below. We finally note that the estimator is in fact polynomial-degree-robust for $hp$ finite element spaces where the polynomial degree is moderately large for elements touching the artificial boundary.

We would also argue that requiring that the elements touching the artificial boundary are large is a generalization of the shape-regularity constraint on the mesh. Indeed, for an interior vertex, this constraint ensures that all element sharing the vertex have comparable size. Similarly, requiring that elements sharing the same vertex have comparable polynomial degrees is a standard assumption, see e.g. [20, Eq. (2.2)]. For the case of a vertex lying on the artificial boundary, we can imagine that the elements touching the boundary have ghost neighboring elements on which we impose that $u_h = 0$. These ghost elements may be considered has being very large and with a low polynomial degree, thus requiring the actual elements in the patch to be large with low polynomial degrees as well. The requirement that $\kappa_K h_K/p_K^+$ is large for elements touching the boundary therefore appear natural.

Second, the term $\|\boldsymbol{\nabla} u\|_K$ may also appear problematic at first glance. However, this term is small if the artificial boundary is well-placed, so that it does not impair the quality of the upper bound. Besides, in the context of an adaptive algorithms, this term can be made smaller by including new elements in the mesh to push the artificial boundary further away in case the elements touching $\Gamma^h$ are marked for refinements. Again, we propose an adaptive algorithm with this feature in Section 6.

We further derive similar estimates for the Helmholtz equation for waveguide problems with perfectly matched layers (PML). In this case, as already happens in bounded media [5, 12], the efficiency bound in (1.2) essentially holds unchanged, but the reliability estimate in (1.1) only holds true if the estimator is pre-multiplied by the problem's inf-sup constant, which is often unavailable. Following [5, 12], we can nevertheless show that the bound in (1.1) essentially remains true on sufficiently refined meshes. In contrast to bounded media, however, the duality techniques employed in the proof are more complex here, since we can not afford to globally refine the mesh.

We complement our theoretical findings with numerical examples that illustrate the key reliability and efficiency bounds. In fact, we propose an adaptive algorithm that automatically refine the mesh and extend $\Omega^h$ iteratively. This is done by employing Dörfler's marking with a slight variation of the usual refinement procedure. In practice, this algorithm seem to deliver optimal convergence rates, although this is not analyzed here.

The remainder of this work is organized as follows. Section 2 gives the notation and preliminary results we employ throughout. In Section 3, we establish our a posteriori error estimates for the reaction diffusion equation, and in Section 4, we do the same for the Helmholtz problem. Section 5 provides a priori error estimates for the Helmholtz problem with partially refined meshes compatible with the efficiency of the error estimator, and we present numerical examples in Section 6. We also provide a set of appendices containing technical results. Specifically, Appendix A contains a trace inequality, we recall the equivalence between the PML Helmholtz problem and its original formulation in Appendix B, and we establish approximation results that are key to our duality techniques in Appendix C.

## 2. Preliminary results and notation

Throughout the manuscript, we employ boldface notation for vector fields. Unless specified, functions and fields in Section 3 are real-valued, whereas they are complex-valued in Section 4.

2.1. **Domain.** We consider a connected open set $\Omega \subset \mathbb{R}^d$, $d = 2$ or 3, that is not necessarily bounded. We assume that the boundary $\Gamma := \partial\Omega$ of $\Omega$ is Lipschitz and that all its connected components correspond to the boundaries of (possibly unbounded) polytopes. Notice that the case $\Gamma = \emptyset$ is not excluded and we do not assume that $\Gamma$ is bounded.

2.2. **Mesh.** We consider a triangulation $\mathcal{T}_h$ of (a subset of) $\Omega$. Specifically, $\mathcal{T}_h$ is a collection of (open) non-degenerate simplices $K$ which are matching in the usual sence (see e.g. [13, Definition 6.1.1]), i.e., if the intersection $\overline{K_-} \cap \overline{K_+}$ of two distinct elements $K_\pm \in \mathcal{T}_h$ is not a empty, then it is a full sub-simplex (i.e. vertex, edge or face) of both elements.

We denote by $\Omega^h$ the open set covered by $\mathcal{T}_h$, i.e., the interior of $\cup_{K\in\mathcal{T}_h}\overline{K}$, and we assume that $\Omega^h$ is connected with a Lischiptz boundary. We denote by $\Gamma^h := \partial\Omega^h \setminus \Gamma$ the artificial boundary introduced by discretization.

For an element $K \in \mathcal{T}_h$, $h_K$ and $\rho_K$ respectively denote the diameters of $K$ and of the largest ball contained in $\overline{K}$. We then denote by $\beta_K := h_K/\rho_K \geq 1$ the shape-regularity parameter of $K$. If $\mathcal{T} \subset \mathcal{T}_h$ is a collection of elements, then $\beta_{\mathcal{T}} := \max_{K\in\mathcal{T}} \beta_K$ is the shape-regularity measure of $\mathcal{T}$.

2.3. **Hidden constants.** If $\mathcal{T} \subset \mathcal{T}_h$ is a set of elements, and $A, B \geq 0$ are two real numbers, then we write that $A \preceq_{\mathcal{T}} B$ when $A \leq C(\beta_{\mathcal{T}})B$, where $C(\beta_{\mathcal{T}})$ is a real number only depending on $\beta_{\mathcal{T}}$. We also use the notation $A \preceq_K B$ when $\mathcal{T} = \{K\}$ for $K \in \mathcal{T}_h$. We also employ the notation $A \lesssim B$ if there exists a constant $C$ independent of $A$, $B$ such that $A \leq CB$. In this case, we will always detail on which other previously introduced symbols the constant $C$ depends, if any.

2.4. **Function spaces.** Throughout this manuscript, for an open set $U \subset \Omega$, $L^2(U)$ is the Lebesgue space of square integrable functions defined on $U$, and we denote by $(\cdot,\cdot)_U$ and $\|\cdot\|_U$ its usual inner product and norm. If $\mathfrak{m} : U \to \mathbb{R}$ is a uniformly positive and bounded measurable function, then the application $L^2(U) \ni v \to \|v\|_{\mathfrak{m},U} := \sqrt{(\mathfrak{m}v, v)_U}$ is a norm on $L^2(U)$ equivalent to $\|\cdot\|_U$. Similarly, $\boldsymbol{L}^2(U) := [L^2(U)]^d$ collects vector-valued function. We use the same notation for the inner-products and norms of $L^2(U)$ and $\boldsymbol{L}^2(U)$. For the weighted $\boldsymbol{L}^2(U)$-norm, the coefficient is a symmetric matrix-valued function, and positivity and boundedness are understood in the sense of quadratic forms.

We will often implicitly extend functions by zero in $\Omega$, meaning that we respectively identify $L^2(U)$ and $\boldsymbol{L}^2(U)$ with subspaces of $L^2(\Omega)$ and $\boldsymbol{L}^2(\Omega)$.

We classically use the notation $H^1(U)$ for the Sobolev space of $v \in L^2(U)$ such that $\boldsymbol{\nabla} v \in \boldsymbol{L}^2(U)$, where $\boldsymbol{\nabla}$ stands for the gradient defined in the sense of distributions. If $\gamma \subset \partial U$ is a relatively open subset of the boundary, then $H^1_\gamma(U)$ contains elements of $H^1(U)$ with vanishing trace on $\gamma$. When $\gamma = \partial U$, we also write $H^1_0(U)$ instead of $H^1_\gamma(U)$.

The key definitions and properties of Lebesgue and Sobolev spaces we shall need in this work are established, e.g., in [22, Chapter 3].

We will also need the vector-valued Sobolev space $\boldsymbol{H}(\mathrm{div}, U)$ of fields $\boldsymbol{v} \in \boldsymbol{L}^2(U)$ such that $\boldsymbol{\nabla} \cdot \boldsymbol{v} \in L^2(U)$, where $\boldsymbol{\nabla}\cdot$ is the divergence operator understood in the sense of distributions. We denote by $\boldsymbol{H}_\gamma(\mathrm{div}, U)$ the elements of $\boldsymbol{H}(\mathrm{div}, U)$ with vanishing normal trace on $\gamma$, i.e., for $\boldsymbol{w} \in \boldsymbol{H}(\mathrm{div}, U)$, we have $\boldsymbol{w} \in \boldsymbol{H}_\gamma(\mathrm{div}, U)$ iff $(\boldsymbol{w}, \boldsymbol{\nabla}\phi)_U + (\boldsymbol{\nabla}\cdot\boldsymbol{w}, \phi)_U = 0$ for all $\phi \in H^1_{\gamma^c}(U)$, with $\gamma_c$ the (open) complement of $\overline{\gamma}$ in $\partial U$. We refer the reader to [16] for more details.

2.5. **Partition of unity on vertex patches.** We denote by $\mathcal{V}_h$ the set of vertices of the mesh $\mathcal{T}_h$. If $K \in \mathcal{T}_h$ is an element, then $\mathcal{V}_h(K) \subset \mathcal{V}_h$ denotes the set of vertices of $K$. If $\boldsymbol{a} \in \mathcal{V}_h$, then $\mathcal{T}_h^{\boldsymbol{a}}$ collects the elements $K \in \mathcal{T}_h$ such that $\boldsymbol{a} \in \mathcal{V}_h(K)$, and $\omega^{\boldsymbol{a}}$ is the open domain covered by $\mathcal{T}_h^{\boldsymbol{a}}$. For $K \in \mathcal{T}_h$, we shall also need $\mathcal{T}_h^K := \cup_{\boldsymbol{a}\in\mathcal{V}_h(K)}\mathcal{T}_h^{\boldsymbol{a}}$, the patch of elements surrounding $K$. We then denote by $\widetilde{K}$ the open domain covered by $\mathcal{T}_h^K$.

For $\boldsymbol{a} \in \mathcal{V}_h$ the "hat function" $\psi^{\boldsymbol{a}}$ is the only element of $\mathcal{P}_1(\mathcal{T}_h) \cap H^1(\Omega^h)$ such that $\psi^{\boldsymbol{a}}(\boldsymbol{b}) = \delta_{\boldsymbol{a},\boldsymbol{b}}$ for all $\boldsymbol{b} \in \mathcal{V}_h$, where $\mathcal{P}_1(\mathcal{T}_h)$ is the set of functions whose restriction to each element $K \in \mathcal{T}_h$ is affine (see Section 2.7 below). Then, $\overline{\omega^{\boldsymbol{a}}}$ is the support of $\psi^{\boldsymbol{a}}$, and the notation $\gamma^{\boldsymbol{a}} := \{\boldsymbol{x} \in \partial\omega^{\boldsymbol{a}} \mid \psi^{\boldsymbol{a}}(\boldsymbol{x}) = 0\}$ and $\gamma^{\boldsymbol{a}}_{\mathrm{c}} := \partial\omega^{\boldsymbol{a}} \setminus \overline{\gamma^{\boldsymbol{a}}}$ will be useful.

In particular, we introduce

$$\boldsymbol{H}_0(\mathrm{div}, \omega^{\boldsymbol{a}}) := \{\boldsymbol{w} \in \boldsymbol{H}(\mathrm{div}, \omega^{\boldsymbol{a}}) \mid \boldsymbol{w}\cdot\boldsymbol{n} = 0 \text{ on } \gamma^{\boldsymbol{a}}\},$$

and

$$H^1_\star(\omega^{\boldsymbol{a}}) := \left\{v \in H^1(\omega^{\boldsymbol{a}}) \mid v = 0 \text{ on } \gamma^{\boldsymbol{a}}_{\mathrm{c}}\right\}.$$

It is important to observe that if $v \in H^1_\star(\omega^{\boldsymbol{a}})$, then $\psi^{\boldsymbol{a}} v \in H^1_0(\omega^{\boldsymbol{a}})$.

We denote by $L^2_\star(\omega^{\boldsymbol{a}}) := \boldsymbol{\nabla}\cdot\boldsymbol{H}_0(\mathrm{div}, \omega^{\boldsymbol{a}})$ the range of the divergence operator. It coincides with $L^2(\omega^{\boldsymbol{a}})$ when $\gamma^{\boldsymbol{a}}_{\mathrm{c}} \neq \emptyset$, and consists of functions with vanishing mean value otherwise.

2.6. **Functional inequalities.** For $K \in \mathcal{T}_h$, $v \in H^1(K)$ and $\nu > 0$, we have the following multiplicative trace inequality:

$$\rho_K^{-1/2}\|v\|_{\partial K} \leq \max\left(\beta_K, \sqrt{d+1}(\nu\rho_K)^{-1}\right)\left(\nu^2\|v\|_K^2 + \|\boldsymbol{\nabla}v\|_K^2\right)^{1/2}. \tag{2.1}$$

This type of inequalities is standard (see e.g. [8, Eq. (4.4)]), but because it is core to the present analysis, we include a proof in Appendix A for completeness.

By combining the product rule and a Poincaré inequality on $H^1_\star(\omega^{\boldsymbol{a}})$, we can show that

$$\|\boldsymbol{\nabla}(\psi^{\boldsymbol{a}}(v - v^{\boldsymbol{a}}))\|_{\omega^{\boldsymbol{a}}} \preceq_{\mathcal{T}_h^{\boldsymbol{a}}} \|\boldsymbol{\nabla}v\|_{\omega^{\boldsymbol{a}}}, \tag{2.2}$$

for all $v \in H^1_\star(\omega^{\boldsymbol{a}})$ where $v^{\boldsymbol{a}} = (v,1)_{\omega^{\boldsymbol{a}}}/(1,1)_{\omega^{\boldsymbol{a}}}$ if $\gamma^{\boldsymbol{a}}_{\mathrm{c}} = \emptyset$ and 0 otherwise. This inequality is established in detail, e.g., in [14, Eq. (3.9)].

2.7. **Polynomial spaces.** For $q \geq 0$ and $K \in \mathcal{T}_h$, $\mathcal{P}_q(K)$ denotes the space of polynomials defined on $K$ and of degree at most $q$. $\boldsymbol{\mathcal{P}}_q(K) := [\mathcal{P}_q(K)]^d$ then collects vector-valued polynomials, and $\boldsymbol{RT}_q(K) := \boldsymbol{\mathcal{P}}_q + \boldsymbol{x}\mathcal{P}_q(K)$ is the Raviart–Thomas polynomial space.

If $\mathcal{T} \subset \mathcal{T}_h$ is a set of mesh elements, then $\mathcal{P}_q(\mathcal{T})$, and $\boldsymbol{RT}_q(\mathcal{T})$ stand for the for the set of functions whose restriction on each $K \in \mathcal{T}$ respectively belongs to $\mathcal{P}_q(K)$ and $\boldsymbol{RT}_q(K)$. Note that these "broken" spaces do not have "build-in" compatibility conditions.

We refer the reader to [13, Chapters 5 and 11] for an extensive description of these spaces.

2.8. **Discrete trace inequality.** Consider an element $K \in \mathcal{T}_h$. Then, for all $q \geq 1$,

$$\rho_K^{1/2}\|v_q\|_{\partial K} \preceq_K q\|v_q\|_K \tag{2.3}$$

for all $v_q \in \mathcal{P}_q(K)$, see [26, Theorem 4.76].

2.9. **Discretization space.** The discretization of the PDEs considered in the work are realized with an $hp$ finite element space. Specifically, for each element $K \in \mathcal{T}_h$, we fix a polynomial degree $p_K \geq 1$, and consider the standard finite element space

$$V_h := \left\{ v_h \in H_0^1(\Omega^h) \mid v_h|_K \in \mathcal{P}_{p_K}(K) \ \forall K \in \mathcal{T}_h \right\}.$$

For each vertex $\boldsymbol{a} \in \mathcal{V}_h$, we will also need the notation $p_-^{\boldsymbol{a}} := \min_{K \in \mathcal{T}_h^{\boldsymbol{a}}} p_K$ and $p_+^{\boldsymbol{a}} := \max_{K \in \mathcal{T}_h^{\boldsymbol{a}}} p_K$. We will similarly need the notation $p_K^+ := \max_{\boldsymbol{a} \in \mathcal{V}_h(V)} p_+^{\boldsymbol{a}}$ and $p_K^- := \min_{\boldsymbol{a} \in \mathcal{V}_h(V)} p_-^{\boldsymbol{a}}$ for $K \in \mathcal{T}_h$.

Although this is trivial, we explicitly remark that $V_h \subset H_0^1(\Omega)$ if we consider that functions are extended by zero. In fact, we shall make this implicit extension by zero often hereafter, and consider $V_h$ as a subset of $H_0^1(\Omega)$ rather that $H_0^1(\Omega^h)$.

2.10. **Broken polynomial projection.** For any $q \geq 1$ and $r \in L^2(\Omega)$, we define $\Pi_q r \in L^2(\Omega^h)$ by requiring that $(\Pi_q r)|_K \in \mathcal{P}_q(K)$ and

$$(r - \Pi_q r, s_q)_K = 0$$

for all $s_q \in \mathcal{P}_q(K)$ and $K \in \mathcal{T}_h$. In this case [26, Theorem 4.76], whenever $r \in H^1(K)$, we have

$$\|r - \Pi_q r\|_K \leq c_{\mathrm{P},K} \frac{h_K}{q} \|\boldsymbol{\nabla} r\|_K, \tag{2.4}$$

with a constant $c_{\mathrm{P},K}$ only depending on $\beta_K$.

## 3. Reaction-diffusion equation

3.1. **Setting.** We consider the bilinear form

$$a(\phi, v) := (\kappa^2 \phi, v)_\Omega + (\boldsymbol{\nabla} \phi, \boldsymbol{\nabla} v)_\Omega$$

for $\phi, v \in H_0^1(\Omega)$, where $\kappa : \Omega \to \mathbb{R}_+$ is a measurable function uniformly bounded away from 0 and $+\infty$. For any open set $U \subset \mathbb{R}^d$, we introduce the norm

$$|\!|\!|\phi|\!|\!|_{\kappa,U}^2 := \|\kappa \phi\|_U^2 + \|\boldsymbol{\nabla} \phi\|_U^2$$

over $H^1(U)$.

Throughout this section, we fix a right-hand side $f \in L^2(\Omega)$. Notice that we do not assume that $f$ is compactly supported.

Due to Lax-Milgram lemma (see e.g. [22, Lemma 2.32]), there exists a unique $u \in H_0^1(\Omega)$ such that

$$a(u, v) = (f, v)_\Omega \qquad \forall v \in H_0^1(\Omega).$$

Similarly, there exists a unique $u_h \in V_h$ such that

$$a(u_h, v_h) = (f, v_h) \qquad \forall v_h \in V_h.$$

3.2. **Piecewise constant coefficient.** For simplicity, we assume that $\kappa$ is piecewise constant, and in particular that $\kappa|_K = \kappa_K \in \mathbb{R}$ for each element $K \in \mathcal{T}_h$. This assumption is in fact not really needed to derive error upper bounds, but it is important to derive lower bounds. Places where the assumption can be lifted are commented upon.

For each $\boldsymbol{a} \in \mathcal{V}_h$ the notation $\kappa^{\boldsymbol{a}} := \max_{K \in \mathcal{T}_h^{\boldsymbol{a}}} \kappa_K$ will be useful. Similarly, we let $\kappa_K^+ := \max_{\boldsymbol{a} \in \mathcal{V}_h(K)} \kappa_K$ for all $K \in \mathcal{T}_h$.

We emphasize here that considering a piecewise constant matrix-valued diffusion coefficient in the definition of the bilinear form $a$ could be done without any fundamental difficulty. For clarity of the exposition, we have chosen not to do so.

3.3. **Prager–Synge type estimate.** We start with a Prager–Synge estimate. In its standard version (see e.g. [4, Eq. (2)], [14, Theorem 3.3], [25] and the references therein), the estimate assumes that $\boldsymbol{\sigma} \cdot \boldsymbol{n} = 0$ on the part of $\partial\Omega^h$ where $u - u_h \neq 0$. Here, we cannot afford such a property for the flux on the artificial boundary $\Gamma^h$, motivating a modification of the estimate. Crucially, we draw the reader's attention to the fact that in general, neither $u$ nor $v$ belongs to $H_0^1(\Omega^h)$ in the theorem below.

**Theorem 3.1** (Prager–Synge estimate)**.** *Let $\boldsymbol{\sigma} \in \boldsymbol{H}(\mathrm{div}, \Omega^h)$ with $(\boldsymbol{\nabla} \cdot \boldsymbol{\sigma}, r_h)_K = (f - \kappa^2 u_h, r_h)_K$ for all $r_h \in \mathcal{P}_{q_K}(K)$, $q_K \geq 1$, for all $K \in \mathcal{T}_h$. Then, we have*

$$
a(u - u_h, v) \leq \sum_{K \in \mathcal{T}_h} \left( c_{\mathrm{P},K} \frac{h_K}{q_K} \|f - \kappa^2 u_h - \boldsymbol{\nabla} \cdot \boldsymbol{\sigma}_h\|_K + \|\boldsymbol{\sigma} + \boldsymbol{\nabla} u_h\|_K \right) \|\boldsymbol{\nabla} v\|_K + |(\boldsymbol{\sigma} \cdot \boldsymbol{n}, v)_{\Gamma^h}| + |(f, v)_{\Omega \setminus \Omega^h}| \tag{3.1}
$$

*for all $v \in H_0^1(\Omega)$.*

Before establishing (3.1), a few remarks are in order.

(i) In principle, the penultimate term in (3.1) should be interpreted as a duality pairing between $H^{-1/2}(\Gamma^h)$ and $H_{00}^{1/2}(\Gamma^h)$. We will, however, only employ this estimate for vector fields with normal traces in $L^2(\Gamma^h)$, so that we keep the inner-product notation.

(ii) Although the sharpest possible constant $c_{\mathrm{P},K}$ for which the upper bound is valid may not be explicitly available, we can always get a guaranteed upper bound by replacing $c_{\mathrm{P},K}$ by $q_K/\pi$. This is done by using the standard Poincaré inequality $\|v - (v,1)_K/(1,1)_K\|_K \leq (h_K/\pi)\|\boldsymbol{\nabla} v\|_K$ (see, e.g., [2]) instead of the Poincaré-like inequality in (2.4).

(iii) In most applications, $f$ typically has compact support and we can immediately design the mesh so that $\operatorname{supp} f \subset \Omega^h$. In this case, the last term in (3.1) may be dropped.

(iv) The estimate in (3.1) holds true under the sole assumption that $\kappa$ is measurable and uniformly bounded away from 0 and $+\infty$. In particular, there is no need to assume that $\kappa$ is piecewise constant here.

*Proof.* We start by writing that

$$
\begin{aligned}
a(u - u_h, v) &= (f, v)_\Omega - \kappa^2 (u_h, v)_{\Omega^h} - (\boldsymbol{\nabla} u_h, \boldsymbol{\nabla} v)_{\Omega^h} \\
&= (f, v)_{\Omega \setminus \Omega^h} + (f - \kappa^2 u_h - \boldsymbol{\nabla} \cdot \boldsymbol{\sigma}, v)_{\Omega^h} + (\boldsymbol{\nabla} \cdot \boldsymbol{\sigma}, v)_{\Omega^h} - (\boldsymbol{\nabla} u_h, \boldsymbol{\nabla} v)_{\Omega^h}.
\end{aligned}
$$

The first term in the right-hand side is already present in (3.1). For the second term, due to our assumptions on $\boldsymbol{\sigma}$, we can write that

$$
\begin{aligned}
(f-\kappa^2 u_h - \boldsymbol{\nabla}\cdot\boldsymbol{\sigma}, v)_{\Omega^h} &= \sum_{K\in\mathcal{T}_h} (f-\kappa^2 u_h - \boldsymbol{\nabla}\cdot\boldsymbol{\sigma}, v - \Pi_{q_K} v)_K \\
&\leq \sum_{K\in\mathcal{T}_h} c_{\mathrm{P},K}\frac{h_K}{q_K}\|f-\kappa^2 u_h - \boldsymbol{\nabla}\cdot\boldsymbol{\sigma}\|_K \|\boldsymbol{\nabla} v\|_K,
\end{aligned}
$$

where we used the Poincaré-like inequality in (2.4). For the remaining terms, we integrate by parts

$$
\begin{aligned}
(\boldsymbol{\nabla}\cdot\boldsymbol{\sigma}, v)_{\Omega^h} - (\boldsymbol{\nabla} u_h, \boldsymbol{\nabla} v)_{\Omega^h} &= (\boldsymbol{\sigma}\cdot\boldsymbol{n}, v)_{\Gamma^h} - (\boldsymbol{\sigma}+\boldsymbol{\nabla} u_h, \boldsymbol{\nabla} v)_{\Omega^h} \\
&\leq \sum_{K\in\mathcal{T}_h} \|\boldsymbol{\sigma}+\boldsymbol{\nabla} u_h\|_K \|\boldsymbol{\nabla} v\|_K + |(\boldsymbol{\sigma}\cdot\boldsymbol{n}, v)_{\Gamma^h}|
\end{aligned}
$$

and (3.1) follows. □

**Corollary 3.2** (Prager–Synge estimate with trace inequalities)**.** *Let $\boldsymbol{\sigma}\in\boldsymbol{H}(\mathrm{div},\Omega^h)$ satisfy the assumption of Theorem 3.1 with $\boldsymbol{\sigma}\cdot\boldsymbol{n}\in L^2(\Gamma^h)$. Then, we have*

$$
\begin{aligned}
a(u-u_h, v) \leq \sum_{K\in\mathcal{T}_h} \left( c_{\mathrm{P},K}\frac{h_K}{q_K}\|f-\kappa^2 u_h - \boldsymbol{\nabla}\cdot\boldsymbol{\sigma}\|_K + \|\boldsymbol{\sigma}+\boldsymbol{\nabla} u_h\|_K + \mu_K \rho_K^{1/2}\|\boldsymbol{\sigma}\cdot\boldsymbol{n}\|_{\partial K\cap\Gamma^h}\right) \|v\|_{\kappa,K} \\
+ \|\kappa^{-1} f\|_{\Omega\setminus\Omega^h} \|v\|_{\kappa,\Omega\setminus\Omega^h}
\end{aligned}
$$

*for all $v\in H^1_0(\Omega)$, where $\mu_K := \max(\beta_K, \sqrt{d+1}(\kappa_K\rho_K)^{-1})$.*

*Proof.* Based on (3.1), we simply need to bound $|(\boldsymbol{\sigma}\cdot\boldsymbol{n}, v)_{\Gamma^h}|$. This is done using the multiplicative trace inequality in (2.1) as follows:

$$
\begin{aligned}
|(\boldsymbol{\sigma}\cdot\boldsymbol{n}, v)_{\Gamma^h}| &= \left| \sum_{K\in\mathcal{T}_h} (\boldsymbol{\sigma}\cdot\boldsymbol{n}, v)_{\partial K\cap\Gamma^h} \right| \leq \sum_{K\in\mathcal{T}_h} \|\boldsymbol{\sigma}\cdot\boldsymbol{n}\|_{\partial K\cap\Gamma_h} \|v\|_{\partial K\cap\Gamma^h} \\
&\leq \sum_{K\in\mathcal{T}_h} \|\boldsymbol{\sigma}\cdot\boldsymbol{n}\|_{\partial K\cap\Gamma_h} \|v\|_{\partial K} = \sum_{K\in\mathcal{T}_h} \rho_K^{1/2}\|\boldsymbol{\sigma}\cdot\boldsymbol{n}\|_{\partial K\cap\Gamma_h} \rho_K^{-1/2}\|v\|_{\partial K} \\
&\leq \sum_{K\in\mathcal{T}_h} \mu_K \rho_K^{1/2}\|\boldsymbol{\sigma}\cdot\boldsymbol{n}\|_{\partial K\cap\Gamma_h} \|v\|_{\kappa,K}.
\end{aligned}
$$

□

**Remark 3.3** (Sharper bounds)**.** *The estimate in Corollary 3.2 has the advantage of being fully-explicit. However, at the price of introducing possibly unknown constant, it may be improved in two ways. (i) First, instead of an element-wise estimate, we could employ a trace inequality from $H^1(\Omega^h)$ to $L^2(\Gamma^h)$, which would save half a negative power of h in the final estimate. (ii) The estimate of Corollary 3.2 is in fact only used for $v = u-u_h$, in which case, we have $\|v\|_{\Gamma^h} = \|u\|_{\Gamma^h}$ since $u_h = 0$ on $\Gamma^h$. If $\Gamma^h$ is sufficiently far away from the support of $f$, the trace of $u$ is expected to be very small, and can sometimes be explicitly estimated.*

3.4. **Construction of the estimator.** For all $\boldsymbol{a} \in \mathcal{V}_h$, we let

$$\boldsymbol{\sigma}_h^{\boldsymbol{a}} := \arg \min_{\substack{\boldsymbol{w}_h \in \boldsymbol{H}_0(\mathrm{div},\omega^{\boldsymbol{a}}) \cap \boldsymbol{RT}_{p_+^{\boldsymbol{a}}+2}(\mathcal{T}_h^{\boldsymbol{a}}) \\ \boldsymbol{\nabla}\cdot\boldsymbol{w}_h = \Pi_{p_+^{\boldsymbol{a}}+2}(\psi^{\boldsymbol{a}} f) - \kappa^2 \psi^{\boldsymbol{a}} u_h - \boldsymbol{\nabla}\psi^{\boldsymbol{a}}\cdot\boldsymbol{\nabla} u_h}} \|\boldsymbol{w}_h + \psi^{\boldsymbol{a}} \boldsymbol{\nabla} u_h\|_{\omega^{\boldsymbol{a}}}. \tag{3.2}$$

We also set

$$\boldsymbol{\sigma}_h := \sum_{\boldsymbol{a}\in\mathcal{V}_h} \boldsymbol{\sigma}_h^{\boldsymbol{a}} \tag{3.3}$$

and

$$f_h := \sum_{\boldsymbol{a}\in\mathcal{V}_h} \Pi_{p_+^{\boldsymbol{a}}+2}(\psi^{\boldsymbol{a}} f). \tag{3.4}$$

For each $K \in \mathcal{T}_h$, our estimator is then defined as

$$\eta_K := c_{\mathrm{P},K} \frac{h_K}{p_K+2} \|f - f_h\|_K + \|\boldsymbol{\sigma}_h + \boldsymbol{\nabla} u_h\|_K + \mu_K \rho_K^{1/2} \|\boldsymbol{\sigma}_h \cdot \boldsymbol{n}\|_{\partial K \cap \Gamma^h}. \tag{3.5}$$

**Proposition 3.4** (Equilibrated flux). *For each vertex $\boldsymbol{a}$ the local problem in (3.2) is well-posed. In addition, the flux defined in (3.3) is equilibrated, i.e., we have*

$$(f - \kappa^2 u_h - \boldsymbol{\nabla}\cdot\boldsymbol{\sigma}_h, s_h)_K = 0 \tag{3.6}$$

*for all $K \in \mathcal{T}_h$ and $s_h \in \mathcal{P}_{p_K+2}(K)$.*

We notice that, when $\kappa$ is not piecewise constant (3.6) remains valid if $\boldsymbol{\sigma}_h^{\boldsymbol{a}}$ is constructed with the divergence constraint $\boldsymbol{\nabla}\cdot\boldsymbol{\sigma}_h^{\boldsymbol{a}} = \Pi_{p_+^{\boldsymbol{a}}+2}(\psi^{\boldsymbol{a}} f - \kappa^2 \psi^{\boldsymbol{a}} u_h - \boldsymbol{\nabla}\psi^{\boldsymbol{a}}\cdot\boldsymbol{\nabla} u_h)$ (compare with (3.2)).

*Proof.* The fact that the local problems are well-posed is standard: When $\gamma^{\boldsymbol{a}} = \partial\omega^{\boldsymbol{a}}$, the compatibility condition is satisfied by Galerkin orthogonality since $\psi^{\boldsymbol{a}} \in \mathcal{V}_h$, see e.g. [14, Section 3.1.3]. To establish (3.6), we first observe that since the $\psi^{\boldsymbol{a}}$ form a partition of unity, we have

$$\boldsymbol{\nabla}\cdot\boldsymbol{\sigma}_h = \sum_{\boldsymbol{a}\in\mathcal{V}_h} \Pi_{p_+^{\boldsymbol{a}}+2}(\psi^{\boldsymbol{a}} f) + \kappa^2 u_h.$$

Then, we fix $K \in \mathcal{T}_h$ and $s_h \in \mathcal{P}_{p_K+2}(K)$ and note that $s_h \in \mathcal{P}_{p_+^{\boldsymbol{a}}+2}(K)$ for all $\boldsymbol{a} \in \mathcal{V}_h(K)$. Due to the locality of the partition of unity, this leads to

$$(f + \kappa^2 u_h - \boldsymbol{\nabla}\cdot\boldsymbol{\sigma}_h, s_h)_K = (f - f_h, s_h)_K = \sum_{\boldsymbol{a}\in\mathcal{V}_h} (\psi^{\boldsymbol{a}} f - \Pi_{p_+^{\boldsymbol{a}}+2}(\psi^{\boldsymbol{a}} f), s_h)_K = 0.$$

□

3.5. **Reliability.** We are now ready to state our main reliability result, which is a direct consequence of Theorem 3.1.

**Theorem 3.5** (Reliability). *For the error estimators $\eta_K$ defined in (3.5), we have*

$$|\!|\!| u - u_h |\!|\!|^2 \leq \sum_{K\in\mathcal{T}_h} \eta_K^2 + \|\kappa^{-1} f\|_{\Omega\setminus\Omega^h}^2. \tag{3.7}$$

It is worth pointing out that the reliability estimate in (3.7) still holds true for generic (i.e. not piecewise constant) coefficients $\kappa$ if the coefficient $\mu_K$ introduced in Corollary 3.2 is modified by taking the minimal value of $\kappa$ in $K$ instead of $\kappa_K$.

3.6. **Efficiency.** Here, we establish lower error bounds. For the "usual" part of the estimator, the standard proof applies without modification. For completeness, it is sketched below. We first establish a "vertex-based" efficiency estimate for each component $\boldsymbol{\sigma}_h^{\boldsymbol{a}}$ of the discrete flux $\boldsymbol{\sigma}_h$.

**Lemma 3.6** (Standard efficiency on patches)**.** *For all $\boldsymbol{a} \in \mathcal{V}_h$, we have*

$$\|\boldsymbol{\sigma}_h^{\boldsymbol{a}} + \psi^{\boldsymbol{a}} \boldsymbol{\nabla} u_h\|_{\omega^{\boldsymbol{a}}} \preceq_{\mathcal{T}_h^{\boldsymbol{a}}} \left(1 + \frac{\kappa^{\boldsymbol{a}} h^{\boldsymbol{a}}}{p_-^{\boldsymbol{a}}}\right) \|u - u_h\|_{\kappa,\omega^{\boldsymbol{a}}}. \tag{3.8}$$

*Proof.* Here, we will need a quasi-interpolation operator $Q_h^{\boldsymbol{a}} : H_0^1(\omega^{\boldsymbol{a}}) \to H_0^1(\omega^{\boldsymbol{a}}) \cap V_h$ such that

$$\|v - Q_h^{\boldsymbol{a}} v\|_{\omega^{\boldsymbol{a}}} \preceq_{\mathcal{T}_h^{\boldsymbol{a}}} \frac{h^{\boldsymbol{a}}}{p_-^{\boldsymbol{a}}} \|\boldsymbol{\nabla} v\|_{\omega^{\boldsymbol{a}}}, \qquad \|\boldsymbol{\nabla}(Q_h^{\boldsymbol{a}} v)\|_{\omega^{\boldsymbol{a}}} \preceq_{\mathcal{T}_h^{\boldsymbol{a}}} \|\boldsymbol{\nabla} v\|_{\omega^{\boldsymbol{a}}}, \tag{3.9}$$

for all $v \in H_0^1(\omega^{\boldsymbol{a}})$, where the notation $H_0^1(\omega^{\boldsymbol{a}}) \cap V_h$ stands for the set of functions in $V_h$ supported in $\omega^{\boldsymbol{a}}$. Such an operator $Q_h^{\boldsymbol{a}}$ is provided, e.g., by [23, Theorem 3.3] when $d = 2$, and by [18, Theorem 4.22] for $d = 3$.

With the operator $Q_h^{\boldsymbol{a}}$ and hand, we start by stating that

$$\|\boldsymbol{\sigma}_h^{\boldsymbol{a}} + \psi^{\boldsymbol{a}} \boldsymbol{\nabla} u_h\|_{\omega^{\boldsymbol{a}}} \preceq_{\mathcal{T}_h^{\boldsymbol{a}}} \min_{\substack{\boldsymbol{\sigma}^{\boldsymbol{a}} \in \boldsymbol{H}_0(\mathrm{div},\omega^{\boldsymbol{a}}) \\ \boldsymbol{\nabla}\cdot\boldsymbol{\sigma}^{\boldsymbol{a}} = \Pi_{p_+^{\boldsymbol{a}}+2}(\psi^{\boldsymbol{a}} f) - \kappa^2 \psi^{\boldsymbol{a}} u_h - \boldsymbol{\nabla}\psi^{\boldsymbol{a}}\cdot\boldsymbol{\nabla} u_h}} \|\boldsymbol{\sigma}^{\boldsymbol{a}} + \psi^{\boldsymbol{a}} \boldsymbol{\nabla} u_h\|_{\omega^{\boldsymbol{a}}}.$$

The proof may be found, e.g., in [4, 7, 14, 15]. The Euler-Lagrange equation for the right-hand side consists in finding $\boldsymbol{\sigma}_\star^{\boldsymbol{a}} \in \boldsymbol{H}_0(\mathrm{div},\omega^{\boldsymbol{a}})$ and $r_\star^{\boldsymbol{a}} \in L_\star^2(\omega^{\boldsymbol{a}})$ such that

$$\left\{\begin{array}{rcl} (\boldsymbol{\sigma}_\star^{\boldsymbol{a}}, \boldsymbol{v})_{\omega^{\boldsymbol{a}}} + (r_\star^{\boldsymbol{a}}, \boldsymbol{\nabla}\cdot\boldsymbol{v})_{\omega^{\boldsymbol{a}}} &=& -(\psi^{\boldsymbol{a}} \boldsymbol{\nabla} u_h, \boldsymbol{v})_{\omega^{\boldsymbol{a}}} \\ (\boldsymbol{\nabla}\cdot\boldsymbol{\sigma}_\star^{\boldsymbol{a}}, q)_{\omega^{\boldsymbol{a}}} &=& (\Pi_{p_+^{\boldsymbol{a}}+2}(\psi^{\boldsymbol{a}} f) - \kappa^2 \psi^{\boldsymbol{a}} u_h - \boldsymbol{\nabla}\psi^{\boldsymbol{a}}\cdot\boldsymbol{\nabla} u_h, q)_{\omega^{\boldsymbol{a}}} \end{array}\right.$$

for all $\boldsymbol{v} \in \boldsymbol{H}_0(\mathrm{div},\omega^{\boldsymbol{a}})$ and $q \in L_\star^2(\omega^{\boldsymbol{a}})$. The first equation readily shows that $r_\star^{\boldsymbol{a}} \in H_\star^1(\omega^{\boldsymbol{a}})$, and integration by parts between $\boldsymbol{H}_0(\mathrm{div},\omega^{\boldsymbol{a}})$ and $H_\star^1(\omega^{\boldsymbol{a}})$ actually gives that

$$\boldsymbol{\sigma}_\star^{\boldsymbol{a}} + \psi^{\boldsymbol{a}} \boldsymbol{\nabla} u_h = \boldsymbol{\nabla} r_\star^{\boldsymbol{a}}.$$

Using the second equation, we further see that

$$\begin{aligned}(\boldsymbol{\nabla} r_\star^{\boldsymbol{a}}, \boldsymbol{\nabla} q)_{\omega^{\boldsymbol{a}}} &= (\Pi_{p_+^{\boldsymbol{a}}+2}(\psi^{\boldsymbol{a}} f), q)_{\omega^{\boldsymbol{a}}} - a(u_h, \psi^{\boldsymbol{a}} q) \\ &= -(\psi^{\boldsymbol{a}} f - \Pi_{p_+^{\boldsymbol{a}}+2}(\psi^{\boldsymbol{a}} f), q)_{\omega^{\boldsymbol{a}}} + a(u - u_h, \psi^{\boldsymbol{a}} q)\end{aligned}$$

for all $q \in H_\star^1(\omega^{\boldsymbol{a}})$ since $\psi^{\boldsymbol{a}} q \in H_0^1(\omega^{\boldsymbol{a}})$. The two terms in the right-hand are treated separately.

(i) First, since $\Pi_{p_+^{\boldsymbol{a}}+2}$ is an orthogonal projection, we have

$$\begin{aligned}-(\psi^{\boldsymbol{a}} f - \Pi_{p_+^{\boldsymbol{a}}+2}(\psi^{\boldsymbol{a}} f), q)_{\omega^{\boldsymbol{a}}} &= -(\psi^{\boldsymbol{a}} f - \Pi_{p_+^{\boldsymbol{a}}+2}(\psi^{\boldsymbol{a}} f), q - \Pi_{p_+^{\boldsymbol{a}}+2} q)_{\omega^{\boldsymbol{a}}} \\ &\preceq_{\mathcal{T}_h^{\boldsymbol{a}}} \frac{h^{\boldsymbol{a}}}{p_+^{\boldsymbol{a}}} \|\psi^{\boldsymbol{a}} f - \Pi_{p_+^{\boldsymbol{a}}+2}(\psi^{\boldsymbol{a}} f)\|_{\omega^{\boldsymbol{a}}} \|\boldsymbol{\nabla} q\|_{\omega^{\boldsymbol{a}}}\end{aligned}$$

where we employed the Poincaré-like inequality in (2.4). On the other hand, we have

$$\begin{aligned}\|\psi^{\boldsymbol{a}} f - \Pi_{p_+^{\boldsymbol{a}}+2}(\psi^{\boldsymbol{a}} f)\|_{\omega^{\boldsymbol{a}}} &= \min_{q_h \in \mathcal{P}_{p_+^{\boldsymbol{a}}+2}(\mathcal{T}_h^{\boldsymbol{a}})} \|\psi^{\boldsymbol{a}} f - q_h\|_{\omega^{\boldsymbol{a}}} \\ &\leq \min_{q_h \in \mathcal{P}_{p_+^{\boldsymbol{a}}+2}(\mathcal{T}_h^{\boldsymbol{a}})} \|\psi^{\boldsymbol{a}} f - \psi^{\boldsymbol{a}} \Pi_{p_+^{\boldsymbol{a}}+1} f\|_{\omega^{\boldsymbol{a}}} \leq \|f - \Pi_{p_+^{\boldsymbol{a}}+1} f\|_{\omega^{\boldsymbol{a}}},\end{aligned}$$

leading to the estimate

$$-(\psi^{\boldsymbol{a}} f - \Pi_{p_+^{\boldsymbol{a}}+2}(\psi^{\boldsymbol{a}} f), q)_{\omega^{\boldsymbol{a}}} \preceq_{\mathcal{T}_h^{\boldsymbol{a}}} \frac{h^{\boldsymbol{a}}}{p_+^{\boldsymbol{a}}} \|f - \Pi_{p_+^{\boldsymbol{a}}+1} f\|_{\omega^{\boldsymbol{a}}} \|\boldsymbol{\nabla} q\|_{\omega^{\boldsymbol{a}}}.$$

(ii) For the remaining term, we first employ Galerkin orthogonality to write that

$$a(u - u_h, \psi^{\boldsymbol{a}} q) = a(u - u_h, w^{\boldsymbol{a}})$$

where $w^{\boldsymbol{a}} = \psi^{\boldsymbol{a}}(q - q^{\boldsymbol{a}})$ with $q^{\boldsymbol{a}} = (q,1)_{\omega^{\boldsymbol{a}}}/(1,1)_{\omega^{\boldsymbol{a}}}$ if $\gamma_{\rm c}^{\boldsymbol{a}} = \emptyset$ and 0 otherwise. Then, by continuity of $a$, we further write

$$a(u - u_h, \psi^{\boldsymbol{a}} q) = a(u - u_h, w^{\boldsymbol{a}} - Q_h^{\boldsymbol{a}} w^{\boldsymbol{a}}) \leq |\!|\!| u - u_h |\!|\!|_{\kappa,\omega^{\boldsymbol{a}}} |\!|\!| w^{\boldsymbol{a}} - Q_h^{\boldsymbol{a}} w^{\boldsymbol{a}} |\!|\!|_{\kappa,\omega^{\boldsymbol{a}}}$$

where we have used the fact that $w^{\boldsymbol{a}} \in H_0^1(\omega^{\boldsymbol{a}})$. Using the estimates in (3.9) and then the inequality in (2.2) shows that

$$|\!|\!| w^{\boldsymbol{a}} - Q_h^{\boldsymbol{a}} w^{\boldsymbol{a}} |\!|\!|_{\kappa,\omega^{\boldsymbol{a}}} \preceq_{\mathcal{T}_h^{\boldsymbol{a}}} \left(1 + \frac{\kappa^{\boldsymbol{a}} h^{\boldsymbol{a}}}{p_-^{\boldsymbol{a}}}\right) \|\boldsymbol{\nabla}(\psi^{\boldsymbol{a}}(q - q^{\boldsymbol{a}}))\|_{\omega^{\boldsymbol{a}}} \preceq_{\mathcal{T}_h^{\boldsymbol{a}}} \left(1 + \frac{\kappa^{\boldsymbol{a}} h^{\boldsymbol{a}}}{p_-^{\boldsymbol{a}}}\right) \|\boldsymbol{\nabla} q\|_{\omega^{\boldsymbol{a}}},$$

leading to (3.8). □

Classicall, the next step in the proof is to recombined the above "vertex-based" estimates to obtain "element-based" lower bounds for the total flux $\boldsymbol{\sigma}_h$.

**Corollary 3.7** (Standard efficiency on elements)**.** *We have*

$$\|\boldsymbol{\sigma}_h + \boldsymbol{\nabla} u_h\|_K \preceq_{\mathcal{T}_h^K} \left(1 + \frac{\kappa_K^+ h_K}{p_K^-}\right) |\!|\!| u - u_h |\!|\!|_{\kappa,\widetilde{K}} \tag{3.10}$$

*for all $K \in \mathcal{T}_h$.*

*Proof.* It simply follows from the triangular inequality

$$\|\boldsymbol{\sigma}_h + \boldsymbol{\nabla} u_h\|_K = \left\| \sum_{\boldsymbol{a} \in \mathcal{V}_h(K)} (\boldsymbol{\sigma}_h^{\boldsymbol{a}} + \psi^{\boldsymbol{a}} \boldsymbol{\nabla} u_h) \right\|_K \leq \sum_{\boldsymbol{a} \in \mathcal{V}_h(K)} \|\boldsymbol{\sigma}_h^{\boldsymbol{a}} + \psi^{\boldsymbol{a}} \boldsymbol{\nabla} u_h\|_{\omega^{\boldsymbol{a}}},$$

and (3.8), where we used that $\max_{\boldsymbol{a} \in \mathcal{V}_h(K)} h^{\boldsymbol{a}} \preceq_{\mathcal{T}_h^K} h_K$ due to the shape-regularity of the mesh. □

We now provide an upper bound for the extra term in the estimator due to the presence of the artificial boundary.

**Lemma 3.8** (Efficiency for the boundary term)**.** *We have*

$$\rho_K^{1/2} \|\boldsymbol{\sigma}_h \cdot \boldsymbol{n}\|_{\partial K \cap \partial \Gamma^h} \preceq_K p_K^+ \left\{ \left(1 + \frac{\kappa_K^+ h_K}{p_K^-}\right) |\!|\!| u - u_h |\!|\!|_{\kappa,K} + \|\boldsymbol{\nabla} u\|_K \right\} \tag{3.11}$$

*for all $K \in \mathcal{T}_h$.*

*Proof.* We start with the discrete trace inequality in (2.3) giving

$$\rho_K^{1/2} \|\boldsymbol{\sigma}_h \cdot \boldsymbol{n}\|_{\partial K \cap \partial \Gamma^h} \preceq_K (p_K^+ + 2) \|\boldsymbol{\sigma}_h\|_K \preceq_K p_K^+ \|\boldsymbol{\sigma}_h\|_K.$$

The estimate in (3.11) then follows from the triangular inequalities

$$\|\boldsymbol{\sigma}_h\|_K \leq \|\boldsymbol{\sigma}_h + \boldsymbol{\nabla} u_h\|_K + \|\boldsymbol{\nabla} u_h\|_K \leq \|\boldsymbol{\sigma}_h + \boldsymbol{\nabla} u_h\|_K + \|\boldsymbol{\nabla}(u - u_h)\|_K + \|\boldsymbol{\nabla} u\|_K$$

and the standard efficiency result in (3.8). □

At that point, the efficiency of the estimator follows by collecting the previous estimates.

**Theorem 3.9** (Efficiency)**.** *For all $K \in \mathcal{T}_h$, we have*

$$\eta_K \preceq_{\mathcal{T}_h^K} \left(1+\delta_{\partial K\cap\Gamma^h}\frac{p_K^+}{\kappa_K h_K}\right)\left(1+\frac{\kappa_K^+ h_K}{p_K^-}\right)\|u-u_h\|_{\kappa,\widetilde{K}} + \frac{h_K}{p_K}\|f-\Pi_{p_K+2}f\|_K + \delta_{\partial K\cap\Gamma^h}\frac{p_K^+}{\kappa_K h_K}\|\boldsymbol{\nabla} u\|_K, \tag{3.12}$$

*where $\delta_{\partial K\cap\Gamma^h} := 1$ if $\partial K\cap\Gamma^h \neq \emptyset$ and $0$ otherwise.*

*Proof.* After observing that $\mu_K \preceq_K (\kappa_K h_K)^{-1}$, we simply combine (3.10) and (3.11). □

We make the following remarks on the efficiency estimate:

(i) Compared to standard efficiency estimates, our error lower bounds have the extra $\|\boldsymbol{\nabla} u\|_K$ term in the right-hand side. This term is however expected to be small if the artificial boundary $\Gamma^h$ is sufficiently far away from the origin. In addition, in an adaptive context, this term can be made smaller by adding new elements to the mesh to push the artificial boundary away from the origin if the estimators associated with elements touching $\Gamma^h$ are large. The numerical examples in Section 6 below seem to indicate that this is a viable strategy.

(ii) The estimator is $p$-robust except at the boundary, where only the local polynomial degree appears. In particular, the proposed estimator is fully robust if $p_K^+ \lesssim 1$ and $\kappa_K h_K \gtrsim 1$ for boundary elements, which is practically doable in an $hp$-adaptive algorithm. In particular, an $h$-adaptive algorithm satisfying this property is presented in the numerical examples in Section 6 below.

(iii) The multiplicative factor $1+\kappa_K^+ h_K/p_K^-$ in (3.12) indicates that the estimator might not be robust for coarse meshes in reaction-dominated regimes. This is standard, and is not due to the fact that unbounded media are considered here. In bounded domains, this limitation can be lifted by modifying the definition of the estimator, see [27, 28]. We believe that such modifications are also possible here.

## 4. Helmholtz problem

4.1. **Settings.** In this section, we describe the Helmholtz problem with PML that we are solving. The motivation of the PML is given in Appendix B.

4.1.1. *Domain.* We assume that the domain $\Omega$ is split into two non-overlapping open sets $\Omega_0, C_\infty \subset \mathbb{R}^d$ such that $\overline{\Omega} = \overline{\Omega_0} \cap \overline{C_\infty}$, where $\Omega_0$ is bounded, and $C_\infty$ consists of a finite number $R \geq 1$ of semi-infinite cylinders. Specifically, we assume that there exist unit vectors $\{\boldsymbol{t}^r\}_{r=1}^R \subset \mathbb{R}^d$, real numbers $\{\ell^r\}_{r=1}^R \subset \mathbb{R}_+$ and (relatively open) Lipschitz connected cross sections $\{\Sigma^r\}_{r=1}^R \subset \mathbb{R}^{d-1}$ such that

$$C_\infty := \bigcup_{r=1}^R C_\infty^r, \qquad C_\infty^r := \left\{\boldsymbol{x}\in\mathbb{R}^d \quad | \quad \boldsymbol{x}\cdot\boldsymbol{t}^r > \ell^r, \quad \boldsymbol{x}-(\boldsymbol{x}\cdot\boldsymbol{t}^r)\boldsymbol{t}^r \in \Sigma^r\right\},$$

and we further assume that the cylinders $C_\infty^r$ are disjoint. The setting is illustrated on Figure 4.1.

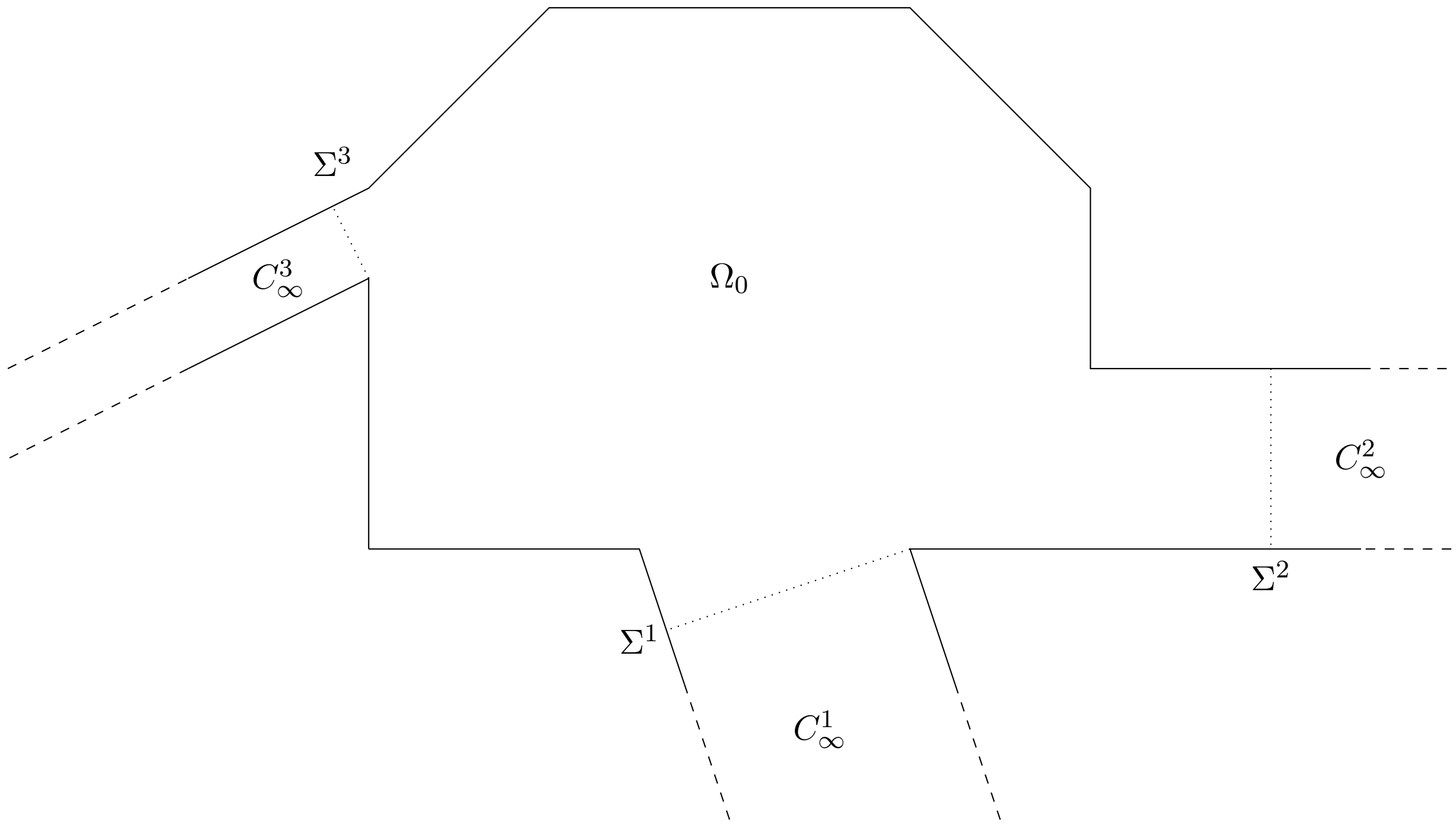

FIGURE 4.1. Example of a waveguide configuration with three cylinders

4.1.2. *Tangential coordinates.* It will be convenient to introduce the vector field $\boldsymbol{t} : \Omega \to \mathbb{R}^d$ by setting $\boldsymbol{t}(\boldsymbol{x}) = \boldsymbol{t}^r$ if $\boldsymbol{x} \in C^r_\infty$ and $\boldsymbol{t}(\boldsymbol{x}) := \boldsymbol{e}_1$ otherwise. The choice we make inside $\Omega_0$ is arbitrary, and it is only introduced to homogenize notation later on. For $\boldsymbol{x} \in \Omega$, then introduce $x_\times = \boldsymbol{x} \cdot \boldsymbol{t}$ and $\boldsymbol{x}_\perp = \boldsymbol{x} - (\boldsymbol{x} \cdot \boldsymbol{t})\boldsymbol{t}$. We will also need the differential operators $\partial_\times v = \boldsymbol{\nabla} v \cdot \boldsymbol{t}$ and $\boldsymbol{\nabla}_\perp v = \boldsymbol{\nabla} v - (\partial_\times v)\boldsymbol{t}$ for $v \in H^1_0(\Omega)$.

4.1.3. *Cut-off frequencies.* For each $r \in \{1, \dots, R\}$, there exists a basis simultaneously orthogonal in $L^2(\Sigma^r)$ and $H^1_0(\Sigma^r)$ given by the eigen modes of the Laplace operator in $\Sigma^r$. Specifically, for $j \in \mathbb{N}$ there exists $\psi^r_j \in H^1_0(\Sigma^r)$ and $\lambda^r_j \geq 0$ such that

$$\left\{ \begin{array}{rcll} -\Delta_\perp \psi^r_j & = & (\lambda^r_j)^2 \psi^r_j & \text{in } \Sigma^r, \\ \psi^r_j & = & 0 & \text{on } \partial\Sigma^r, \end{array} \right. \tag{4.1}$$

and the family $(\psi^r_j)_{j \geq 0}$ is an orthonormal basis of $L^2(\Sigma^r)$, and an orthogonal basis $H^1_0(\Sigma^r)$.

Throughout the remainder of this work, we assume that $k$ is not a cut-off frequency, which means that $k \neq \lambda^r_j$ for all $r \in \{1, \dots, R\}$ and $j \in \mathbb{N}$. Equivalently, we have

$$k_\star := \min_{1 \leq r \leq R} \min_{j \geq 0} |k^r_j| > 0,$$

and we introduce the notation

$$\mu := \frac{k_\star}{k} > 0.$$

4.1.4. *Damping coefficient.* We consider a fixed, user-defined damping $\gamma \in \mathbb{C}$ such that $\operatorname{Re} \gamma \geq 1$ and $\operatorname{Im} \gamma \geq 1$. We introduce the shorthand notation $\gamma_{\rm r} := \operatorname{Re} \gamma$ and $\gamma_{\rm i} := \operatorname{Im} \gamma$, as well as $\gamma_\star := \min(\gamma_{\rm r}, \gamma_{\rm i}) \geq 1$.

It will be useful to introduce the coefficients $\alpha := \mathbf{1}_{\Omega_0} + \gamma \mathbf{1}_{C_\infty}$, and $\underline{\boldsymbol{A}} := \alpha^{-1} \boldsymbol{t} \otimes \boldsymbol{t} + \alpha(\underline{\boldsymbol{I}} - \boldsymbol{t} \otimes \boldsymbol{t})$, so that $\alpha = 1$ and $\underline{\boldsymbol{A}} = \underline{\boldsymbol{I}}$ on $\Omega_0$. for all $\phi, v \in H^1_0(\Omega)$. We will frequently use the short-hand otation $\underline{\boldsymbol{A}}_{\rm r} := \operatorname{Re} \underline{\boldsymbol{A}}$.

4.1.5. *Helmholtz problem.* We introduce the sesquilinear form

$$\begin{aligned} b(\phi, v) &:= -k^2 (\alpha \phi, v)_\Omega + (\alpha \boldsymbol{\nabla}_\perp \phi, \boldsymbol{\nabla}_\perp v)_\Omega + (\alpha^{-1} \partial_\times \phi, \partial_\times v)_\Omega \\ &= -k^2 (\alpha \phi, v)_\Omega + (\underline{\boldsymbol{A}} \boldsymbol{\nabla} \phi, \boldsymbol{\nabla} v)_\Omega \end{aligned}$$

for all $\phi, v \in H^1_0(\Omega)$.

Our model problem then consists in finding $u \in H^1_0(\Omega)$ such that

$$b(u, v) = (f, v)_\Omega \tag{4.2}$$

for all $v \in H^1_0(\Omega)$, where $f \in L^2(\Omega_0)$ is a fixed right-hand side supported in $\Omega_0$.

4.1.6. *Energy norm.* For an open set $U \subset \Omega$ and $\phi \in H^1(U)$, we consider the weighted norm

$$|||\phi|||^2_{k,U} := k^2 \|\phi\|^2_{\alpha_{\rm r}, U} + \|\boldsymbol{\nabla} \phi\|^2_{\underline{\boldsymbol{A}}_{\rm r}, U}.$$

which allows us the write the following Gårding inequality

$$\operatorname{Re} b|_U(\phi, \phi) = -k^2 \|\phi\|^2_{\alpha_{\rm r}, U} + \|\boldsymbol{\nabla} \phi\|^2_{\underline{\boldsymbol{A}}_{\rm r}, U} = |||\phi|||^2_{k,U} - 2k^2 \|\phi\|^2_{\alpha_{\rm r}, U}, \tag{4.3}$$

where $b|_U$ stands for the sesquilinear form defined as $b$ but with integrals taken over $U$ instead of $\Omega$.

We also observe that

$$\|\phi\|_{\alpha_{\rm r}^{-1}, U} \leq \|\phi\|_{\alpha_{\rm r}, U} \qquad k^2 \|\phi\|^2_U + \|\boldsymbol{\nabla} \phi\|^2_U \leq \frac{|\gamma|^2}{\gamma_{\rm r}} |||\phi|||^2_{k,U}. \tag{4.4}$$

for all $\phi \in L^2(U)$.

4.1.7. *Well-posedness.* We will assume that the problem in (4.2) is well-posed, and specifically, that there exists a real number $\beta_{\rm st} > 0$ such that

$$\min_{\substack{\phi \in H^1_0(\Omega) \\ |||\phi|||_{k,\Omega} = 1}} \max_{\substack{\phi \in H^1_0(\Omega) \\ |||\phi|||_{k,\Omega} = 1}} \operatorname{Re} b(\phi, v) = \frac{1}{\beta_{\rm st}} \tag{4.5}$$

As we detail in Appendix B, this holds true under natural assumptions. Namely, that the original Helmholtz problem (where a radiation condition is set at infinity rather than the PML) is uniquely solvable.

4.1.8. *Discrete solution.* In the remainder of this section we consider an arbitrary element $u_h \in V_h$ such that $b(u_h, v_h) = (f, v_h)$ for all $v_h \in V_h$. Note that we even though we do not know a priori that such a discrete solution exist, it is a satisfactory assumption to make in a posteriori analysis where we assumed that the problem has been successfully solved anyway. We will discuss the existence of discrete solution later on in Section 5 below.

4.2. **Prager–Synge estimate.** We are now ready to establish a Prager–Synge estimate on the PDE residual. This is done analogously to the coercive case analyzed in Section 3.

**Theorem 4.1** (Prager–Synge estimate)**.** *Let* $\boldsymbol{\sigma} \in \boldsymbol{H}(\mathrm{div}, \Omega^h)$ *with* $(\boldsymbol{\nabla} \cdot \boldsymbol{\sigma}, r_h)_K = (f + k^2 \alpha u_h, r_h)_K$ *for all* $r_h \in \mathcal{P}_{q_K}(K)$, $q_K \geq 1$, *for all* $K \in \mathcal{T}_h$. *Assume further that* $\boldsymbol{\nabla} \cdot \boldsymbol{\sigma} = k^2 \alpha u_h$ *in* $C_\infty$. *Then, for all* $v \in H^1_0(\Omega)$, *we have*

$$|b(u-u_h, v)| \leq \sum_{K \in \mathcal{T}_h} \left( c_{\mathrm{P},K} \frac{h_K}{q_K} \|f + k^2 \alpha u_h - \boldsymbol{\nabla} \cdot \boldsymbol{\sigma}\|_K + \|\boldsymbol{\sigma} + \underline{\boldsymbol{A}} \boldsymbol{\nabla} u_h\|_{\underline{\boldsymbol{A}}_{\mathrm{r}}^{-1},K} \right) \|\boldsymbol{\nabla} v\|_{\underline{\boldsymbol{A}}_{\mathrm{r}},K} + |(\boldsymbol{\sigma} \cdot \boldsymbol{n}, v)_{\Gamma^h}|. \tag{4.6}$$

*Proof.* We start with Stokes' integration by part formula:

$$\begin{aligned} b(u-u_h, v) &= (f, v)_\Omega + k^2 (\alpha u_h, v)_\Omega - (\underline{\boldsymbol{A}} \boldsymbol{\nabla} u_h, \boldsymbol{\nabla} v)_\Omega \\ &= (f + k^2 \alpha u_h - \boldsymbol{\nabla} \cdot \boldsymbol{\sigma}, v)_\Omega - (\boldsymbol{\sigma} + \underline{\boldsymbol{A}} \boldsymbol{\nabla} u_h, \boldsymbol{\nabla} v)_\Omega + (\boldsymbol{\sigma} \cdot \boldsymbol{n}, v)_{\partial\Omega}. \end{aligned}$$

The first term is treated as in Theorem 3.1. Introducing $f_h := \boldsymbol{\nabla} \cdot \boldsymbol{\sigma} - k^2 \alpha u_h$, we have

$$\begin{aligned} |(f + k^2 \alpha u_h - \boldsymbol{\nabla} \cdot \boldsymbol{\sigma}, v)_{\Omega^h}| &\leq \sum_{K \in \mathcal{T}_h} c_{\mathrm{P},K} \frac{h_K}{q_K} \|f - f_h\|_K \|\boldsymbol{\nabla} v\|_K \\ &= \sum_{K \in \mathcal{T}_h} c_{\mathrm{P},K} \frac{h_K}{q_K} \|f - f_h\|_K \|\boldsymbol{\nabla} v\|_{\underline{\boldsymbol{A}}_{\mathrm{r}},K} \end{aligned}$$

since the terms where $\underline{\boldsymbol{A}}|_K \neq \underline{\boldsymbol{I}}$ vanish by assumption, as $\operatorname{supp} f \subset \Omega_0$. For the second term, we simply have

$$|(\boldsymbol{\sigma} + \underline{\boldsymbol{A}} \boldsymbol{\nabla} u_h, \boldsymbol{\nabla} v)_\Omega| \leq \sum_{K \in \mathcal{T}_h} |(\boldsymbol{\sigma} + \underline{\boldsymbol{A}} \boldsymbol{\nabla} u_h, \boldsymbol{\nabla} v)_K| \leq \sum_{K \in \mathcal{T}_h} \|\boldsymbol{\sigma} + \underline{\boldsymbol{A}} \boldsymbol{\nabla} u_h\|_{\underline{\boldsymbol{A}}_{\mathrm{r}}^{-1},K} \|\boldsymbol{\nabla} v\|_{\underline{\boldsymbol{A}}_{\mathrm{r}},K},$$

and the last term is easily bounded since $v = 0$ on $\partial\Omega \setminus \Gamma^h$. □

4.3. **Estimator.** The Prager–Synge estimate established above leads to a natural definition of the error estimator. Again, we follow the lines of Section 3.

For all $\boldsymbol{a} \in \mathcal{V}_h$, we let

$$\boldsymbol{\sigma}_h^{\boldsymbol{a}} := \arg \min_{\substack{\boldsymbol{w}_h \in \boldsymbol{H}_0(\mathrm{div}, \omega^{\boldsymbol{a}}) \cap \boldsymbol{RT}_{p_+^{\boldsymbol{a}}+2}(\mathcal{T}_h^{\boldsymbol{a}}) \\ \boldsymbol{\nabla} \cdot \boldsymbol{w}_h = \Pi_{p_+^{\boldsymbol{a}}+2}(\psi^{\boldsymbol{a}} f) + k^2 \psi^{\boldsymbol{a}} \alpha u_h - \underline{\boldsymbol{A}} \boldsymbol{\nabla} \psi^{\boldsymbol{a}} \cdot \boldsymbol{\nabla} u_h}} \|\boldsymbol{w}_h + \psi^{\boldsymbol{a}} \underline{\boldsymbol{A}} \boldsymbol{\nabla} u_h\|_{\underline{\boldsymbol{A}}_{\mathrm{r}}^{-1}, \omega^{\boldsymbol{a}}}. \tag{4.7}$$

We also set

$$\boldsymbol{\sigma}_h := \sum_{\boldsymbol{a} \in \mathcal{V}_h} \boldsymbol{\sigma}_h^{\boldsymbol{a}} \tag{4.8}$$

and

$$f_h := \sum_{\boldsymbol{a} \in \mathcal{V}_h} \Pi_{p_+^{\boldsymbol{a}}+2}(\psi^{\boldsymbol{a}} f). \tag{4.9}$$

For each $K \in \mathcal{T}_h$, our estimator is then defined as

$$\eta_K := c_{\mathrm{P},K} \frac{h_K}{p_K + 2} \|f - f_h\|_K + \|\boldsymbol{\sigma}_h + \underline{\boldsymbol{A}} \boldsymbol{\nabla} u_h\|_{\underline{\boldsymbol{A}}_{\mathrm{r}}^{-1},K} + \mu_K \rho_K^{1/2} \|\boldsymbol{\sigma}_h \cdot \boldsymbol{n}\|_{\partial K \cap \Gamma^h}. \tag{4.10}$$

where

$$\mu_K := \frac{|\gamma|^2}{\gamma_{\mathrm{r}}} \max(\beta_K, \sqrt{d+1}(k\rho_K)^{-1}). \tag{4.11}$$

**Proposition 4.2** (Properties of the estimator)**.** *The local problem in (3.2) is well-posed for all $\boldsymbol{a} \in \mathcal{V}_h$, and the the flux defined in (3.3) is equilibrated, i.e., we have*

$$(f + k^2\alpha u_h - \boldsymbol{\nabla}\cdot\boldsymbol{\sigma}_h, q_h)_K = 0 \tag{4.12}$$

*for all $K \in \mathcal{T}_h$ and $q_h \in \mathcal{P}_{p_K+2}(K)$. The resulting estimator controls the residual, i.e., for all $v \in H^1_0(\Omega)$, we have*

$$|b(u-u_h, v)| \leq \eta \|v\|_{k,\Omega}. \tag{4.13}$$

*Proof.* The proof of (4.12) very closely follows to the one of Proposition 3.4 in the coercive case. For (4.13), recalling (4.4), we follow the proof of Corollary 3.2 and simply use the alternative estimate

$$\begin{aligned}
\rho_K^{1/2}\|\boldsymbol{\sigma}\cdot\boldsymbol{n}\|_{\partial K}\rho_K^{-1/2}\|v\|_{\partial K} &\leq \max(\beta_K, \sqrt{d+1}(k\rho_K)^{-1})\rho_K^{1/2}\|\boldsymbol{\sigma}\cdot\boldsymbol{n}\|_{\partial K}\left(k^2\|v\|_K^2 + \|\boldsymbol{\nabla} v\|_K^2\right)^{1/2} \\
&\leq \frac{|\gamma|^2}{\gamma_{\mathrm{r}}}\max(\beta_K, \sqrt{d+1}(k\rho_K)^{-1})\rho_K^{1/2}\|\boldsymbol{\sigma}\cdot\boldsymbol{n}\|_{\partial K}\|v\|_{k,K} \\
&= \mu_K\rho_K^{1/2}\|\boldsymbol{\sigma}\cdot\boldsymbol{n}\|_{\partial K}\|v\|_{k,K}
\end{aligned}$$

for all $K \in \mathcal{T}_h$ and $v \in H^1(K)$ at the very end of the proof. □

4.4. **Basic reliability and efficiency.** Given Proposition 4.2, we can immediately provide basic reliability and efficiency results.

**Theorem 4.3** (Coarse error estimate)**.** *We have*

$$\|u-u_h\|_{k,\Omega} \leq \beta_{\mathrm{st}}\eta. \tag{4.14}$$

*Proof.* Since $u - u_h \in H^1_0(\Omega)$, there exists $v^\star \in H^1_0(\Omega)$ with $\|v^\star\|_{k,\Omega} = 1$ such that

$$\|u-u_h\|_{k,\Omega} \leq \beta_{\mathrm{st}} \operatorname{Re} b(u-u_h, v^\star),$$

and (4.14) then follows from (4.13). □

**Theorem 4.4** (Efficiency)**.** *For all $K \in \mathcal{T}_h$, we have*

$$\begin{aligned}
\eta_K \lesssim \left(1 + \delta_{\partial K\cap\Gamma_h}\frac{p_K^+}{kh_K}\right)\left(1 + \frac{kh_K}{p_K^-}\right)\|u-u_h\|_{k,\widetilde{K}} \\
+ \frac{h_K}{p_K}\|f - \Pi_{p_K+2}f\|_K + \delta_{\partial K\cap\Gamma^h}\frac{p_K^+}{kh_K}\|\boldsymbol{\nabla} u\|_K,
\end{aligned} \tag{4.15}$$

*where the hidden constant depends on $\gamma$ and $\beta_{\mathcal{T}_h^K}$.*

The estimate in (4.15) is proved by following the lines of Theorem 3.9, see also [5, Lemma 3.7]. For shortness, we do not reproduce the proof. The dependency on $\gamma$ in (4.15) is due to the coefficient $\mu_K$ introduced in (4.11) and the constants appears in the norm equivalence between $\|\cdot\|_{\underline{\boldsymbol{A}}_{\mathrm{r}},U}$ and the standard $\boldsymbol{L}^2(U)$ norm. For brevity, we do not track the dependency, but it is expected to be mild (polynomial in $|\gamma|$ and $1/\gamma_\star$ with small powers).

The reliability estimate in (4.14) is not fully satisfactory. Indeed, it involves the inf-sup stability constant $\beta_{\mathrm{st}}$ which is frequency-dependant and generally unknown. In the next two sections, we show that this estimate can be refined if the mesh suitably refined.

4.5. **Approximation factor.** Following [5, 12], we introduce an "approximation factor" to measure the refinement of the mesh. We will need some extra notation.

4.5.1. *Notation.* For each $r \in \{1, \dots, R\}$ and $n > 0$, we let

$$C_n^r := \left\{ \boldsymbol{x} \in C_\infty^r \mid x_\times - \ell^r < \frac{n}{k} \right\},$$

i.e., the $n$ first wavelengths are contained in $C_\infty^r$ (note that $n$ does not to be an integer). We then consider the domain $\Omega_n$ such that

$$\overline{\Omega_n} := \overline{\Omega_0} \cup \bigcup_{r=1}^{R} \overline{C_n^r},$$

which consists of $\Omega_0$ extended by $n$ wavelengths in each semi-infinite cylinder. We will also need the notation

$$\Sigma_n := \partial\Omega_n \setminus \partial\Omega, \qquad \Omega_n^{\mathrm{c}} := \Omega \setminus \overline{\Omega_n}.$$

For simplicity, we will assume in Theorem 4.5 below that there exists $\nu_{\min} > 0$ such that

$$k \operatorname{diam}(\Sigma^r) \geq \nu_{\min} \tag{4.16}$$

for all $r \in \{1, \dots, R\}$. This assumption is only made to simplify the expression of the upper bound on the approximation factor, and is not mandatory. Indeed, then $\mu \lesssim 1$ where the hidden constant only depends on $\nu_{\min}$.

4.5.2. *Approximation factor.* For $n \geq 0$, we are now ready to introduce the approximation factor

$$\delta_n := \sup_{\substack{\phi \in L^2(\Omega_n) \\ \|\phi\|_{\alpha,\Omega}=1}} \min_{v_h \in V_h} \|\xi_\phi - v_h\|_{k,\Omega} \tag{4.17}$$

where $\xi_\phi$ is the unique element of $H_0^1(\Omega)$ such that

$$b(w, \xi_\phi) = k(\alpha w, \phi)_{\Omega_n}.$$

As we show in Appendix C, the approximation factor becomes small if the mesh is sufficiently refined inside $\Omega_0$ and sufficiently far away in the PML. Crucially, our assumption does not rule out coarse elements far away from $\Omega_0$, which is key to the efficiency of the estimator.

**Theorem 4.5** (Approximation factor)**.** *Assume for simplicity that $k \geq k_0 > 0$. Then, there exists $s > 0$ depending on $\Omega$ and $\gamma$ such that for all $n \geq 0$ and $m \geq n + 1$, we have*

$$\delta_n \lesssim \beta_{\mathrm{st}} \left( (kh_m)^s + kh_m + \frac{e^{-\gamma_\star \mu (m-n)}}{\sqrt{\gamma_\star \mu}} \right) \tag{4.18}$$

*with*

$$h_n := \max_{\substack{K \in \mathcal{T}_h \\ K \cap \Omega_{n+1} \neq \emptyset}} h_K,$$

*where the hidden constant depends on $k_0$, $\Omega$, $\gamma$, $s$, $\nu_{\min}$ and $\beta_{\mathcal{T}_h}$.*

**Remark 4.6** (Comparison with standard duality arguments)**.** *In a standard duality argument, the approximation factor is defined with right-hand sides that are supported in the whole domain $\Omega$ (i.e. $\delta_\infty$ in the present notation), see e.g. [5]. Although for a bounded domain, this is the most straightforward strategy, it requires the mesh to be globally refined everywhere in $\Omega$. Such an approach is therefore not amenable to unbounded domains, since it would* *require infinitely many mesh elements. Using duality techniques with right-hand sides locally supported was already explored in* [17, 19]*, but with slightly different purposes in mind.*

4.6. **Asymptotically constant-free estimates.** Equipped with the definition of the approximation factor, we can proceed with our duality techniques. The first step is an essentially standard Aubin-Nitsche trick, whose proof is included below for completeness.

**Lemma 4.7** (Aubin-Nitsche trick)**.** *For all $n \geq 0$, we have*

$$k\|u - u_h\|_{\alpha,\Omega_n} \leq \delta_n \eta. \tag{4.19}$$

*Proof.* For shortness, we write $e := u - u_h$. Recalling the definition of $\delta_m$ in (4.17), we have

$$b(w, \xi_e) = k(\alpha w, e)_{\Omega_n}$$

for all $w \in H^1_0(\Omega)$. Picking $w = e$ yields

$$k\|e\|^2_{\alpha,\Omega_n} = \operatorname{Re} b(e, \xi_e) = \operatorname{Re} b(e, \xi_e - \Pi_h \xi_e) \leq \eta \|\!|\xi_e - \Pi_h \xi_e|\!\|_{k,\Omega},$$

and (4.19) follows since

$$\|\!|\xi_e - \Pi_h \xi_e|\!\|_{k,\Omega} \leq \delta_m \|e\|_{\alpha,\Omega_n}$$

by definition. □

For our asymptotically constant-free upper bound, we will also need the following cutoff function.

**Lemma 4.8** (Cutoff-function)**.** *Consider the function $\chi_k : \mathbb{R} \to \mathbb{R}$ defined by*

$$\chi_k(x) := \left| \begin{array}{ll} 1 & kx < 0 \\ (1 - kx)^2 & 0 < kx < 1 \\ 0 & 1 < kx. \end{array} \right.$$

*Then, $\chi$ is Lipschitz-continuous, and we have $|\chi_k'| \leq 2k\chi_k^{1/2}$.*

We can now establish the main result of this section.

**Theorem 4.9** (Asymptotically constant-free estimate)**.** *Assume that $\Omega^h \supset \Omega_1$. Then, we have*

$$(1 - 2\delta_1)\|\!|u - u_h|\!\|^2_{k,\Omega_0} \leq (1 + 6\delta_1 + 4\delta_1^2)\eta^2. \tag{4.20}$$

*Proof.* Consider the cutoff function $\chi$ defined by $\chi = 1$ in $\Omega_0$ and $\chi(\boldsymbol{x}) = \chi_k(x_\times - \ell^r)$ in each $C^r$, for $1 \leq r \leq R$. Observe that with this definition, we have $\chi = 1$ on $\Omega_0$, $\chi = 0$ on $\Omega_1^{\mathrm{c}}$, and the following properties hold true

$$\chi \leq \chi^{1/2} \qquad |\partial_\times \chi| \leq 2k\chi^{1/2} \leq 2k. \tag{4.21}$$

in $\Omega$. Besides, we have

$$\|\!|e|\!\|^2_{k,\Omega_0} \leq k^2\|e\|^2_{\chi\alpha_{\mathrm{r}},\Omega} + \|\boldsymbol{\nabla} e\|^2_{\chi\underline{\boldsymbol{A}}_{\mathrm{r}},\Omega} := \operatorname{Re}\left\{k^2(\chi\alpha e, e)_\Omega + (\chi\underline{\boldsymbol{A}}\boldsymbol{\nabla} e, \boldsymbol{\nabla} e)_\Omega\right\}. \tag{4.22}$$

We then expand

$$\begin{aligned} b(e, \chi e) &= -k^2(\alpha e, \chi e)_\Omega + (\alpha \boldsymbol{\nabla}_\perp e, \boldsymbol{\nabla}_\perp(\chi e))_\Omega + (1/\alpha \partial_\times e, \partial_\times(\chi e))_\Omega \\ &= -k^2(\chi\alpha e, e)_\Omega + (\chi\alpha\boldsymbol{\nabla}_\perp e, \boldsymbol{\nabla}_\perp e)_\Omega + (\chi/\alpha\partial_\times e, \partial_\times e)_\Omega + (\partial_\times\chi/\alpha\partial_\times e, e)_\Omega, \end{aligned}$$

where we used the fact that $\boldsymbol{\nabla}_\perp \chi = \boldsymbol{0}$. We find that

$$\operatorname{Re} b(e, \chi e) = -k^2\|e\|^2_{\chi\alpha_{\mathrm{r}},\Omega} + \|\boldsymbol{\nabla} e\|^2_{\chi\underline{\boldsymbol{A}}_{\mathrm{r}},\Omega} + \operatorname{Re}(\partial_\times\chi/\alpha\partial_\times e, e)_\Omega,$$

and

$$k^2\|e\|^2_{\chi\alpha_{\rm r},\Omega} + \|\boldsymbol{\nabla} e\|^2_{\chi\underline{\boldsymbol{A}}_{\rm r},\Omega} = \operatorname{Re} b(e,\chi e) + 2k^2\|e\|^2_{\chi\alpha_{\rm r},\Omega} - \operatorname{Re}(\partial_\times\chi/\alpha\partial_\times e, e)_\Omega.$$

Since $\partial_\times\chi$ is supported in $\Omega_1$, we can control the last term as follows

$$|\operatorname{Re}(\partial_\times\chi/\alpha\partial_\times e, e)_\Omega| = |\operatorname{Re}(\partial_\times\chi/\alpha\partial_\times e, e)_{\Omega_1}| \leq \frac{2k}{|\gamma|}(\chi^{1/2}|\partial_\times e|, |e|)_{\Omega_1}$$
$$\leq 2k\|\partial_\times e\|_{\operatorname{Re}(\chi/\alpha),\Omega}\|e\|_{\alpha_{\rm r},\Omega_1} \leq \delta_1\|\partial_\times e\|^2_{\operatorname{Re}(\chi/\alpha),\Omega} + \delta_1^{-1}\|e\|_{\alpha_{\rm r},\Omega_1},$$

where we used that $\operatorname{Re}(\chi/\alpha) \leq \chi$ as $|\gamma| \geq 1$, (4.21) and the usual algebraic inequality. This leads to

$$(1-\delta_1)\left(k^2\|e\|^2_{\chi\alpha_{\rm r},\Omega} + \|\boldsymbol{\nabla} e\|^2_{\chi\underline{\boldsymbol{A}}_{\rm r},\Omega}\right) \leq |b(e,\chi e)| + (2+\delta_1^{-1})k^2\|e\|^2_{\alpha_{\rm r},\Omega_1}.$$

We now invoke the bound on $\|e\|_{\alpha_{\rm r},\Omega_1}$ in (4.19), leading to

$$(1-\delta_1)\left(k^2\|e\|^2_{\chi\alpha_{\rm r},\Omega} + \|\boldsymbol{\nabla} e\|^2_{\chi\underline{\boldsymbol{A}}_{\rm r},\Omega}\right) \leq |b(e,\chi e)| + (\delta_1 + 2\delta_1^2)\eta^2. \tag{4.23}$$

We now need to focus on the first term in the right-hand side of (4.23). We start by with (4.6), which gives

$$|b(e,\chi e)| \leq \eta\|\boldsymbol{\nabla}(\chi e)\|_{\underline{\boldsymbol{A}}_{\rm r},\Omega} \leq \eta(\|\chi\boldsymbol{\nabla} e\|_{\underline{\boldsymbol{A}}_{\rm r},\Omega} + \|\boldsymbol{\nabla}\chi e\|_{\underline{\boldsymbol{A}}_{\rm r},\Omega}) \tag{4.24}$$

since $\chi e \in H^1_0(\Omega^h)$ due to our assumption on $\Omega^h$. On the one hand, we can write that

$$\|\chi\boldsymbol{\nabla} e\|_{\underline{\boldsymbol{A}}_{\rm r},\Omega} \leq \|\boldsymbol{\nabla} e\|_{\chi\underline{\boldsymbol{A}}_{\rm r},\Omega}$$

since $\chi \leq \chi^{1/2}$. On the other hand, we have

$$\|\boldsymbol{\nabla}\chi e\|_{\underline{\boldsymbol{A}}_{\rm r},\Omega} = \|e\partial_\times\chi\|_{\operatorname{Re}(\alpha^{-1}),\Omega_1} \leq 2k\|e\|_{\operatorname{Re}(\alpha^{-1}),\Omega_1} \leq 2k\|e\|_{\alpha_{\rm r},\Omega_1} \leq 2\delta_1\eta,$$

where we used (4.4). Collecting these last two estimates in (4.24) gives

$$|b(e,\chi e)| \leq \eta\|\boldsymbol{\nabla} e\|_{\chi\underline{\boldsymbol{A}}_{\rm r},\Omega} + 2\delta_1\eta^2. \tag{4.25}$$

We then combine (4.23) with (4.25)

$$(1-\delta_1)\left(k^2\|e\|^2_{\chi\alpha_{\rm r},\Omega} + \|\boldsymbol{\nabla} e\|^2_{\chi\underline{\boldsymbol{A}}_{\rm r},\Omega}\right) \leq \eta\|\boldsymbol{\nabla} e\|_{\chi\underline{\boldsymbol{A}}_{\rm r},\Omega} + (2\delta_1^2 + 3\delta_1)\eta^2,$$

and we employ an algebraic inequality on the first term in the right-hand side

$$(1-2\delta_1)\left(k^2\|e\|^2_{\chi\alpha,\Omega} + \|\boldsymbol{\nabla} e\|^2_{\chi\underline{\boldsymbol{A}}_{\rm r},\Omega}\right) \leq (1+6\delta_1+4\delta_1^2)\eta^2.$$

At that point, (4.20) immediately follows from (4.22). □

This conclude our a posteriori error analysis of the Helmholtz problem. So far, we have assumed the existence of $u_h$, which is a perfectly fine assumption as far as a posteriori error analysis is concerned. Nevertheless, our estimator is only efficient on meshes with coarse elements close to $\Gamma^h$, and the standard a priori theory does not guarantee that a discrete solution exist (not to mention its accuracy) on such meshes. The next section is dedicated to establishing that this is in fact the case, under mesh refinement assumptions that are compatible with the efficiency of the estimator.

## 5. Quasi-optimality on partially refined meshes

The goal of this section is to established the existence and uniqueness of solutions to the discrete Helmholtz problem which consists in finding $u_h \in V_h$ such that

$$b(u_h, v_h) = (f, v_h)_\Omega \tag{5.1}$$

for all $v_h \in V_h$. In fact, we are also going to derive a quasi-optimality estimate, which in turns imply convergence of $u_h$ to $u$ for suitably refined meshes.

Crucially, our a posteriori error analysis requires meshes than are not refined in the vicinity of the artificial boundary $\Gamma^h$. However, standard results based on the Schatz duality argument require globally refined meshes. Hereafter, we therefore propose a modification of the proof that does not hinge on any mesh refinement close to the $\Gamma^h$. The proof does demand, however, that the artificial boundary is placed sufficiently far away from $\Omega_0$.

As already pointed out in Remark 4.6, some of the ideas employed in this section have been previously considered in [17, 19]. For completeness though, we provide a full treatment here.

The starting point of our analysis is to estimate the norm of the solution in the PML by exploiting its damping properties. This results in a Gårding inequality where, crucially, the $L^2$ norm in the right-hand side is not taken over the PML. In what for $q \geq 0$, we introduce the short-hand notation $b_q := b|_{\Omega_q}$.

**Lemma 5.1** (Gårding inequality in the PML)**.** *The estimate*

$$\frac{2\gamma_{\rm r}^2}{|\gamma|^2}\|\partial_\times \phi\|^2_{\alpha_{\rm r}^{-1},C_q} \leq \operatorname{Re} b_q(\phi, -i\frac{\gamma}{\gamma_{\rm i}}\phi) + k^2\|\phi\|^2_{\Omega_0} \tag{5.2}$$

*holds true for all $q \geq 0$ and $\phi \in H^1_0(\Omega_q)$.*

*Proof.* Let $\phi \in H^1_0(\Omega)$. Simple manipulations show that

$$\begin{aligned} b_q(\phi, \gamma\phi) &= -k^2(\alpha\overline{\gamma}\phi, \phi)_{\Omega_q} + (\alpha\overline{\gamma}\boldsymbol{\nabla}_\perp\phi, \boldsymbol{\nabla}_\perp\phi)_{\Omega_q} + (\overline{\gamma}/\alpha\partial_\times\phi, \partial_\times\phi)_{\Omega_q} \\ &= -k^2\overline{\gamma}\|\phi\|^2_{\Omega_0} + \overline{\gamma}\|\boldsymbol{\nabla}\phi\|^2_{\Omega_0} - k^2|\gamma|^2\|\phi\|^2_{C_q} + |\gamma|^2\|\boldsymbol{\nabla}_\perp\phi\|^2_{C_q} + \frac{\overline{\gamma}}{\gamma}\|\partial_\times\phi\|^2_{C_q}, \end{aligned}$$

and therefore

$$\begin{aligned} \operatorname{Im} b_q(\phi, \gamma\phi) &= k^2\gamma_{\rm i}\|\phi\|^2_{\Omega_0} - \gamma_{\rm i}\|\boldsymbol{\nabla}\phi\|^2_{\Omega_0} + \operatorname{Im}\frac{\overline{\gamma}^2}{|\gamma|^2}\|\partial_\times\phi\|^2_{C_\infty} \\ &= k^2\gamma_{\rm i}\|\phi\|^2_{\Omega_0} - \gamma_{\rm i}\|\boldsymbol{\nabla}\phi\|^2_{\Omega_0} - \frac{2\gamma_{\rm i}\gamma_{\rm r}}{|\gamma|^2}\|\partial_\times\phi\|^2_{C_\infty}. \end{aligned}$$

We eventually rewrite this identity as

$$\operatorname{Re} b_q(\phi, -i\frac{\gamma}{\gamma_{\rm i}}\phi) = \frac{2\gamma_{\rm r}^2}{|\gamma|^2}\|\partial_\times\phi\|^2_{1/\gamma_{\rm r},C_\infty} + \|\boldsymbol{\nabla}\phi\|^2_{\Omega_0} - k^2\|\phi\|^2_{\Omega_0}.$$

and (5.2) follows by moving the last term of the right-hand side to the left-hand side. □

We are now ready to established a new Gårding inequality, where the energy norm is controlled by the sesquilnear form, the $L^2(\Omega_0)$ norm, and the $L^2(\Sigma_q)$ norm (which vanishes for arguments in $V_h$ whenever $\Omega^h \subset \Omega_q$). Crucially, the $L^2(C_\infty)$ norm does not appear in the right-hand side of (5.3), which subsequently enables an improved duality argument.

**Theorem 5.2** (Gårding inequality). *For all $q \geq 1$ and $\phi \in H_0^1(\Omega)$, we have*

$$\|\phi\|_{k,\Omega_q}^2 \leq \operatorname{Re} b_q(\phi, \iota_q \phi) + (2 + \rho_q^2) k^2 \|\phi\|_{\Omega_0}^2 + 2\rho_q k \|\phi\|_{\Sigma_q}^2, \tag{5.3}$$

*where*

$$\rho_q := 2|\gamma| q, \quad \iota_q := 1 - i\frac{\gamma}{\gamma_{\mathrm{i}}}\rho_q^2.$$

*Proof.* We start with the following Poincaré inequality

$$k^2 \|\phi\|_{C_q}^2 \leq 2qk\|\phi\|_{\Sigma_q}^2 + 4q^2 \|\partial_\times \phi\|_{C_q}^2,$$

which is established in Lemma C.1 below. We then write that

$$\begin{aligned} k^2 \|\phi\|_{\alpha_{\mathrm{r}},C_q}^2 &\leq 2\gamma_{\mathrm{r}} qk \|\phi\|_{\Sigma_q}^2 + 4\gamma_{\mathrm{r}} q^2 \|\partial_\times \phi\|_{C_q}^2 \\ &\leq 2\gamma_{\mathrm{r}} qk \|\phi\|_{\Sigma_q}^2 + 4\gamma_{\mathrm{r}}^2 q^2 \|\partial_\times \phi\|_{\alpha_{\mathrm{r}}^{-1},C_q}^2. \end{aligned}$$

Combining this estimate with (5.2), we have

$$k^2 \|\phi\|_{\alpha_{\mathrm{r}},C_q}^2 \leq 2\gamma_{\mathrm{r}} qk \|\phi\|_{\Sigma_q}^2 + 4\frac{|\gamma|^2}{2} q^2 \left( \operatorname{Re} b(\phi, -i\frac{\gamma}{\gamma_{\mathrm{i}}}\phi) + k^2 \|\phi\|_{\Omega_0}^2 \right). \tag{5.4}$$

To finish the proof, we rewrite the standard Gårding inequality in (4.3) as follows:

$$\|\phi\|_{k,\Omega_q}^2 = \operatorname{Re} b_q(\phi,\phi) + 2k^2 \|\phi\|_{\alpha_{\mathrm{r}},\Omega_q}^2 = \operatorname{Re} b_q(\phi,\phi) + 2k^2 \|\phi\|_{\Omega_0}^2 + 2k^2 \|\phi\|_{\alpha_{\mathrm{r}},C_q}^2.$$

Now, (5.3) follows by controlling the last term with (5.4) and regrouping the terms. □

We now continue our Schatz duality argument using our original Gårding inequality. As usual, we will need the continuity of $b$

$$|b(\phi, v)| \leq \frac{|\gamma|}{\gamma_{\mathrm{r}}} \|\phi\|_{k,\Omega} \|v\|_{k,\Omega} \tag{5.5}$$

as well as the "approximation factor" $\delta_0$ defined in (4.17) which quantifies whether the mesh size is sufficiently small for the finite element solution to be quasi-optimal. The key to our modified analysis is that the right-hand side in (4.17) is only supported in $\Omega_0$, whereas the standard argument would take the supremum over $L^2(\Omega)$. Crucially, this guarantees that $\delta_0$ can is small on meshes not refined in the vicinity of $\Gamma^h$, provided that $\Gamma^h$ is sufficiently far away from $\Omega_0$, as shown in Theorem 4.5.

**Lemma 5.3** (Aubin-Nitsche trick). *For all discrete solutions $u_h \in V_h$ to (5.1), we have*

$$k\|e\|_{\Omega_0} \leq \frac{|\gamma|}{\gamma_{\mathrm{r}}} \delta_0 \|u - u_h\|_{k,\Omega}. \tag{5.6}$$

*Proof.* For shortness, we set $e := u - u_h$. Besides, in what follows, we denote by $\xi$ the unique lement of $H_0^1(\Omega)$ such that

$$b(w, \xi) = k(w, e)_{\Omega_0} \tag{5.7}$$

for all $w \in H_0^1(\Omega)$. By selecting the test function $w = e$ in (5.7), we obtain

$$k\|e\|_{\Omega_0}^2 = b(e, \xi) = b(e, \xi - \xi_h)$$

where we employed Galerkin orthogonality, and

$$\xi_h := \arg \min_{v_h \in V_h} \|\xi - v_h\|_{k,\Omega}.$$

Recalling the continuity property of $b$ in (5.5) and the definition of $\delta_0$ in (4.17), we then obtain that

$$k\|e\|_{\Omega_0}^2 \leq \frac{|\gamma|}{\gamma_{\rm r}}\|e\|_{k,\Omega}\|\xi - \xi_h\|_{k,\Omega} \leq \frac{|\gamma|}{\gamma_{\rm r}}\delta_0\|e\|_{k,\Omega}\|e\|_{\Omega_0}, \leq \frac{1}{2k}\frac{|\gamma|^2}{\gamma_{\rm r}^2}\delta_0^2\|e\|_{k,\Omega}^2 + \frac{k}{2}\|e\|_{\Omega_0}^2,$$

and we obtain

$$k^2\|e\|_{\Omega_0}^2 \leq \frac{|\gamma|^2}{\gamma_{\rm r}^2}\delta_0^2\|e\|_{k,\Omega}^2$$

by shifting the last term to the left-hand side and the inequality by $2k$. □

We are now ready to conclude our analysis and establish a quasi-optimality result provided that $\delta_0$ is sufficiently small.

**Theorem 5.4** (Asymptotic quasi-optimality)**.** *Assume that $\Omega^h \subset \Omega_q$ for some $q \geq 1$ and that*

$$(2+\rho_q^2)\frac{|\gamma|^2}{\gamma_{\rm r}^2}\delta_0^2 < 1.$$

*Then, there exists a unique discrete solution $u_h \in V_h$ to (5.1). In addition, there holds*

$$\left(1-(2+\rho_q^2)\frac{|\gamma|^2}{\gamma_{\rm r}^2}\delta_0^2\right)^{1/2}\|u-u_h\|_{k,\Omega} \leq M_q \min_{v_h \in V_h}\|u-v_h\|_{k,\Omega} \leq \frac{M_q\delta_0}{k}\|f\|_{\Omega_0} \tag{5.8}$$

*with*

$$M_q^2 := 15\frac{|\gamma|^2}{\gamma_{\rm r}}|\iota_q|.$$

**Remark 5.5** (Smallnest assumption)**.** *We emphasize that the smallnesst constraint $\rho_q\delta_0 \lesssim 1$ can ideed be met. For that, we assume that for some $m \geq 0$ sufficiently large, $kh_m$ is sufficiently small. We also want to ensure that $kh_K \sim 1$, while ensuring that $\Omega^h \subset \Omega_q$ for $q$ sufficiently small. Let us first observe that the condition $\rho_q\delta_0 \leq \tau < 1$ may be summarized as*

$$q\left((kh_m)^s + kh_m + e^{-\gamma_\star\mu m}\right) \lesssim \tau$$

*where the hidden constant does not depend on $h, m, q$ nor $\tau$. Assuming for simplcity that $kh_m \leq \varepsilon \leq 1$ and picking $m = |\ln(\varepsilon)|/(\mu\gamma_\star)$ we can further simplify the expression as*

$$q\varepsilon^s \lesssim \tau.$$

*We finally need to understand how small $q$ may be taken. For element $K \in \mathcal{T}_h$ in the region $x_\times = k(m+1)$ we have $kh_K \sim \varepsilon$. By using geometric coarsening, we can construct the mesh in such a way that $kh_K \sim 1$ for elements $K \in \mathcal{T}_h$ in the region $x_\times = k(m+r)$ for a fixed $r \geq 2$. Since $m$ grows with $\varepsilon$, we can therefore assume that $q = 2m$ for $\varepsilon$ sufficiently small, and we end up with the condition*

$$|\ln\varepsilon|\varepsilon^s \lesssim \tau,$$

*which can indeed be met for all $\tau > 0$ by selecting $\varepsilon$ small enough.*

*Proof.* Let us denote by

$$v_h := \arg\min_{w_h \in V_h}\|u-w_h\|_{k,\Omega}$$

the best approximation of $u$. For the moment, let $u_h \in V_h$ be any discrete solution. As above, we write $e := u-u_h$ for shortness and $E := u - v_h$. By combining (5.3) with $\phi = e$, and (5.6), we have

$$\|e\|_{k,\Omega_q}^2 \leq \operatorname{Re} b_q(e,\iota_q e) + 2\rho_q k\|e\|_{\Sigma_q}^2 + (2+\rho_q^2)\frac{|\gamma|^2}{\gamma_{\rm r}^2}\delta_0^2\left(\|e\|_{k,\Omega_q}^2 + \|u\|_{k,\Omega_q^{\rm c}}^2\right).$$

where we used the fact that $u_h = 0$ on $\Omega_q^{\mathrm{c}}$. Working on the first term, we write that

$$\operatorname{Re} b_q(e, \iota_q e) = \operatorname{Re} b(e, \iota_q e) + b|_{\Omega_q^{\mathrm{c}}}(e, \iota_q e) = \operatorname{Re} b(e, \iota_q E) + b|_{\Omega_q^{\mathrm{c}}}(u, \iota_q u) \leq \frac{|\gamma|}{\gamma_{\mathrm{r}}} |\iota_q| \left( \|e\|_{k,\Omega} \|E\|_{k,\Omega} + \|u\|^2_{\Omega_q^{\mathrm{c}}} \right).$$

where we employed Galerkin orthogonality and the continuty of $b$ stated in (5.5). Since

$$(2 + \rho_q^2) \frac{|\gamma|^2}{\gamma_{\mathrm{r}}^2} \delta_0^2 \leq 1, \qquad \frac{|\gamma|}{\gamma_{\mathrm{r}}} |\iota_q| \geq 1,$$

by assumption, this leads to

$$\|e\|^2_{k,\Omega_q} \leq \frac{|\gamma|}{\gamma_{\mathrm{r}}} |\iota_q| \|e\|_{k,\Omega} \|E\|_{k,\Omega} + 2\rho_q k \|u\|^2_{\Sigma_q} + (2 + \rho_q^2) \frac{|\gamma|^2}{\gamma_{\mathrm{r}}^2} \delta_0^2 \|e\|^2_{k,\Omega_q} + 2 \frac{|\gamma|}{\gamma_{\mathrm{r}}} |\iota_q| \|u\|^2_{k,\Omega_q^{\mathrm{c}}},$$

By adding $\|u\|^2_{k,\Omega_q^{\mathrm{c}}}$ to both sides, we have

$$\|e\|^2_{k,\Omega} \leq \frac{|\gamma|}{\gamma_{\mathrm{r}}} |\iota_q| \|e\|_{k,\Omega} \|E\|_{k,\Omega} + 2\rho_q k \|u\|^2_{\Sigma_q} + (2 + \rho_q^2) \frac{|\gamma|^2}{\gamma_{\mathrm{r}}^2} \delta_0^2 \|e\|^2_{k,\Omega_q} + 3 \frac{|\gamma|}{\gamma_{\mathrm{r}}} |\iota_q| \|u\|^2_{k,\Omega_q^{\mathrm{c}}},$$

and therefore

$$\left\{ 1 - (2 + \rho_q^2) \frac{|\gamma|^2}{\gamma_{\mathrm{r}}^2} \delta_0^2 \right\} \|e\|^2_{k,\Omega} \leq \frac{|\gamma|}{\gamma_{\mathrm{r}}} |\iota_q| \|e\|_{k,\Omega} \|E\|_{k,\Omega} + 2\rho_q k \|u\|^2_{\Sigma_q} + 3 \frac{|\gamma|}{\gamma_{\mathrm{r}}} |\iota_q| \|u\|^2_{k,\Omega_q^{\mathrm{c}}}.$$

By using an algebraic inequality, we further obtain the estimate

$$\left\{ 1 - (4 + 2\rho^2) \frac{|\gamma|^2}{\gamma_{\mathrm{r}}^2} \delta_0^2 \right\} \|e\|^2_{k,\Omega} \leq \frac{|\gamma|^2}{\gamma_{\mathrm{r}}^2} |\iota_q|^2 \|E\|^2_{k,\Omega} + 4\rho_q k \|u\|^2_{\Sigma_q} + 6 \frac{|\gamma|}{\gamma_{\mathrm{r}}} |\iota_q| \|u\|^2_{k,\Omega_q^{\mathrm{c}}}.$$

We conclude the proof with the trace inequality in (C.2), i.e.,

$$4\rho_q k \|u\|^2_{\Sigma_q} = 4\rho_q k \|E\|^2_{\Sigma_q} \leq 8\rho_q \frac{|\gamma|^2}{\gamma_{\mathrm{r}}} \|E\|^2_{k,\Omega}$$

and the fact that

$$\|u\|^2_{k,\Omega_q^{\mathrm{c}}} \leq \|E\|^2_{k,\Omega},$$

since $\operatorname{supp} v_h \subset \Omega_q$. □

## 6. Numerical examples

In this section, we illustrate the main theoretical findings with numerical examples.

6.1. **Setting.** Throughout we consider two-dimensional domains $\Omega \subset \mathbb{R}^2$ and finite element spaces with a spatially fixed polynomial degree $p$.

6.1.1. *Meshes.* Our meshes are based on infinite Cartesian grids with mesh size 1 align so that the origin is the vertex of a square. Each square is then subdivided into four triangles through its barycenter to produce a conforming simplicial mesh. We then consider refinements of this infinite triangular mesh by newest vertex bisection [3] (NVB). So far, these meshes have an infinite number of elements. We thus introduce an integer parameter $L > 0$ and only select the elements in these meshes that are inside $\Omega$ and such that $\max_{\boldsymbol{x} \in K} |\boldsymbol{x}|_\infty \leq L$, where $|\boldsymbol{x}|_\infty = \max(|x_1|, |x_2|)$. Notice that this truncation never cuts elements (and similarly, we have designed the domains $\Omega$ in such a way that they are nicely covered). Examples of meshes we employ are depicted on Figures 6.2 and 6.4 below.

6.1.2. *Adaptive loop.* We employ our error estimator within an adaptive loop, following the usual

$$\textbf{SOLVE} \implies \textbf{ESTIMATE} \implies \textbf{MARK} \implies \textbf{REFINE}$$

procedure. The modules **SOLVE**, **ESTIMATE** and **MARK** are standard. Specfcially, (a) given a mesh $\mathcal{T}_h$ we compute the associated discrete solution $u_h$ with the `mumps` direct solver [1]. (b) We compute our a posteriori error estimator as per (3.2) and (4.7). (c) We mark elements for refinement using the standard Dörfler marking [11], i.e., we order the elements in descending order according to the estimator and collect the first $\mathcal{M}_h \subset \mathcal{T}_h$ elements in a way that

$$\sum_{K\in\mathcal{M}_h} \eta_K^2 \geq \theta \sum_{K\in\mathcal{T}_h} \eta_K^2,$$

where $\theta := 0.2$.

Our **REFINE** module differs from the standard version. Specifically, if an element $K \in \mathcal{T}_h$ not touching $\Gamma^h$ is marked for refinement, we classically employ the NVB [3] algorithm to refine it. However, when an element touching the boundary is marked for refinement, we do not refine it, and instead increase $L$ by 1 for the next mesh.

We note that this algorithm should ensure that elements touching the artificial remain unrefined, except for possible propagation of refinement of interior elements to ensure conformity in the NVB. In practice, in the examples below, we observe that the elements touching the artificial boundary indeed remain unrefined.

In all the examples below, we start with a subset of the unit Cartesian grid discribed in Section 6.1.1 above. The grid is naturally restricted to $\Omega$ when necassary (for the Helmholtz problem), and the truncation parameter $L$ is specificied for each test case below. We then run the adaptive loop for 64 iterations.

6.1.3. *Error measurements.* The analytical solution of the examples below are not readily available. When considering a discrete solution $u_h$ computed on a mesh $\mathcal{T}_h$ with polynomial degree $p$ and truncation parameter $L$, we approximate the error as

$$|||u - u_h||| \sim |||\widetilde{u} - u_h|||,$$

where $\widetilde{u}_h$ is computed on the same mesh $\mathcal{T}_h$ with polynomial degree $p+2$ and a large, fixed, truncation parameter $\widetilde{L} \gg L$.

6.2. **Reaction diffusion equation.** Here, we first consider the reaction diffusion equation on $\mathbb{R}^2$ with $\kappa := 1$. The right-hand side is the set function $f := \mathbf{1}_{(-1,1)^2}$ and the initial mesh exactly covers its support, i.e., it is aligned so that the origin is a vertex and $L = 1$. We consider the polynomial degrees $p = 1$ and 3.

We note that due to the corner in support of $f$, the solution is $H^2$ regular, but not $H^4$ regular, so that uniform refinements lead to optimal convergence when $p = 1$, but not for $p = 3$. The sequence of meshes generated by the adaptive algorithm for $p = 1$ and 3 are respectively depicted on Figures 6.1 and 6.2. We note that in the 5 (resp. 10) first iterations for $p = 1$) (resp. 3), the algorithm almost only pushes the boundary without refining any elements. This is to be expected since we start with a very small domain $\Omega^h$, with $\Gamma^h$ corresponding to the boundary of the support of $f$. In later iterations, the algorithm properly interleaves refinements of elements and extension of $\Omega^h$. We see in particular that the refinements are concentrated towards the origin, with gradings at the corners of the support of $f$ for $p = 3$, but not for $p = 1$. This is in line with our intuitive expectations. We also see that the boundary is pushed away less often for $p = 1$ than for $p = 3$, which also matches our intuition.

A more precise description of the performance of the adaptive algorithm and the error estimator is provided on Figure 6.3. On the left panel, we see that the estimator $\eta$ and the error converge at the optimal rate $N_{\rm dofs}^{p/2}$, where $N_{\rm dofs}$ is the dimension of the finite element space $V_h$. We also see that the estimator is indeed a guaranteed upper bound for the error. We have also represented the "standard" error estimator $\widetilde{\eta}$, which does not account for the artificial boundary. We see that, perhaps unsurprisingly, it severely underestimate the error on the initial meshes, where the artificial boundary is placed very close to the support of the right-hand side. The right panel of Figure 6.3 present the effectivity indices of our estimator $\eta$ and the standard estimator $\widetilde{\eta}$. There again, we see that our estimator provides a guaranteed upper bound. We also see that this bound is sharp, since the effectivity index is very close to one starting from iteration 10.

6.3. **Helmholtz problem.** We now consider a waveguide problem modeled by the Helmholtz equation with a PML. The domain is "T-shaped", as can be easily seen on Figure 6.4. Rigorously it defined as

$$\Omega = (-\infty, +\infty) \times (-1, 1) \cup (-1, 1) \times (1, -\infty).$$

We select the damping parameter $\gamma = 1 + i$ in the PML, which we impose for $\boldsymbol{x} \in \Omega$ with $|\boldsymbol{x}|_\infty > 5$, i.e.

$$\Omega_0 = (-5, 5) \times (-1, 1) \cup (-1, 1) \times (1, -5).$$

We consider two different values for the wavenumber, namely $k = 0.7 \cdot 2\pi$ and $1.7 \cdot 2\pi$. The right-hand side models the injection of the second guided mode of the waveguide from the bottom branch. Specifically, we let

$$U(\boldsymbol{x}) := e^{iKx_2} \sin(\pi x_1)$$

with $K = \sqrt{k^2 - \pi^2}$ be the guided mode, and we let

$$f := (-k^2 - \Delta)(\chi U)$$

where $\chi$ is a smooth transition functions such that $\chi = 1$ if $x_2 > -3$ and 0 whenever $x_2 < -3.5$. Notice that then $f$ is smooth and supported on $(-1, 1) \times (-3.5, -3)$.

The initial mesh is correctly aligned with the boundary of $\Omega$ and the interface of the PML. Specifically, we take $L = 7$ so that there is four squares in each connected region of the PML in the initial mesh.

On Figure 6.4, we represent the sequence of meshes generated by the adaptive algorithm for the lower frequency. Here too, these refinements intuitively make sense, since the meshes become graded toward the reentrant corners, and are more refined in the domain and at the beginning of the PML than towards the end.

Figure 6.5 then represents the behaviour of the whole error, the error measured only on $\Omega_0$ and our error estimator. On the left panels, we observe that the optimal convergence rate is achieved. On the right panels, we see that on coarse meshes, the estimator underestimate the error, especially for the larger frequency value. This is inline with our theory. As the meshes gets refined, the estimator asymptotically becomes a guaranteed upper bound for the error measured in $\Omega_0$, in accordance with Theorem 4.9. We also see that the effectivity index for the whole error seem to tends to one, as would be the case for the Helmholtz equation in a bounded domain. Although this not behaviour is not completely surprising, it is not covered by the present analysis. We finally note that the artificial boundary is pushed more often for the lower frequency. This intuitively makes sense, since less refinements are required to resolve the higher wavelength in this case.

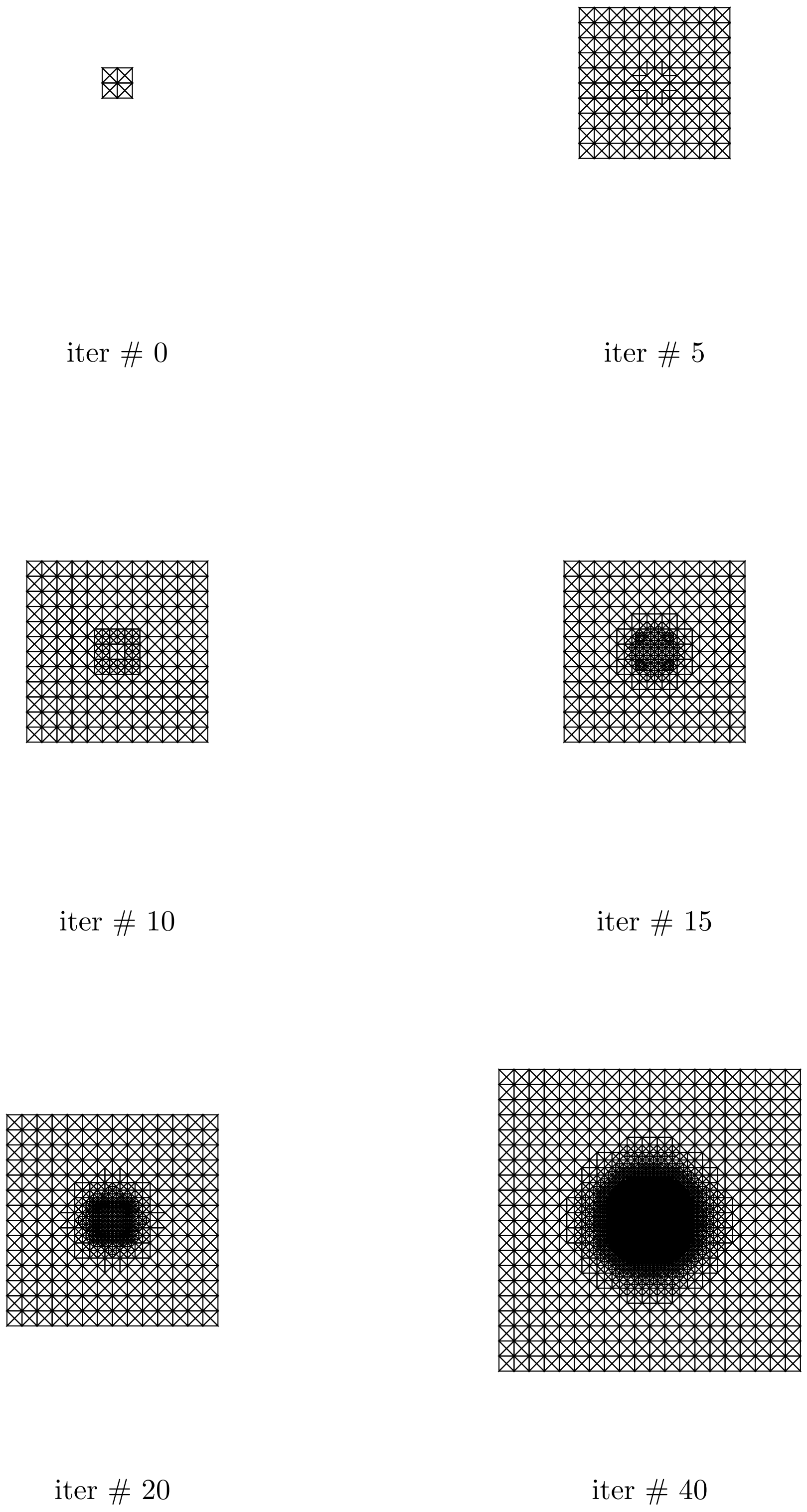

FIGURE 6.1. Sequence of adapted meshes for the reaction diffusion problem with $p = 1$.

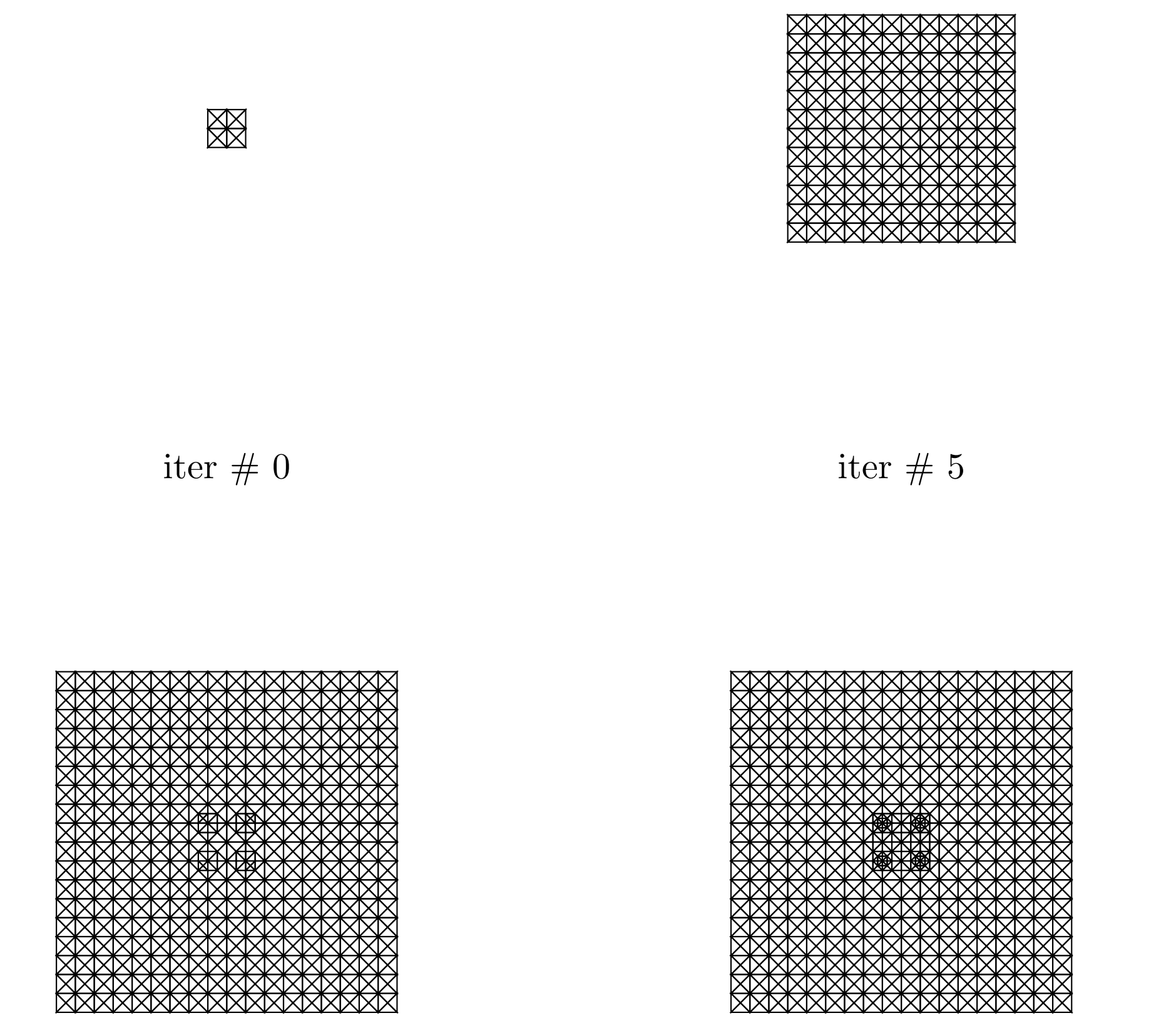

FIGURE 6.2. Sequence of adapted meshes for the reaction diffusion problem with $p = 3$.

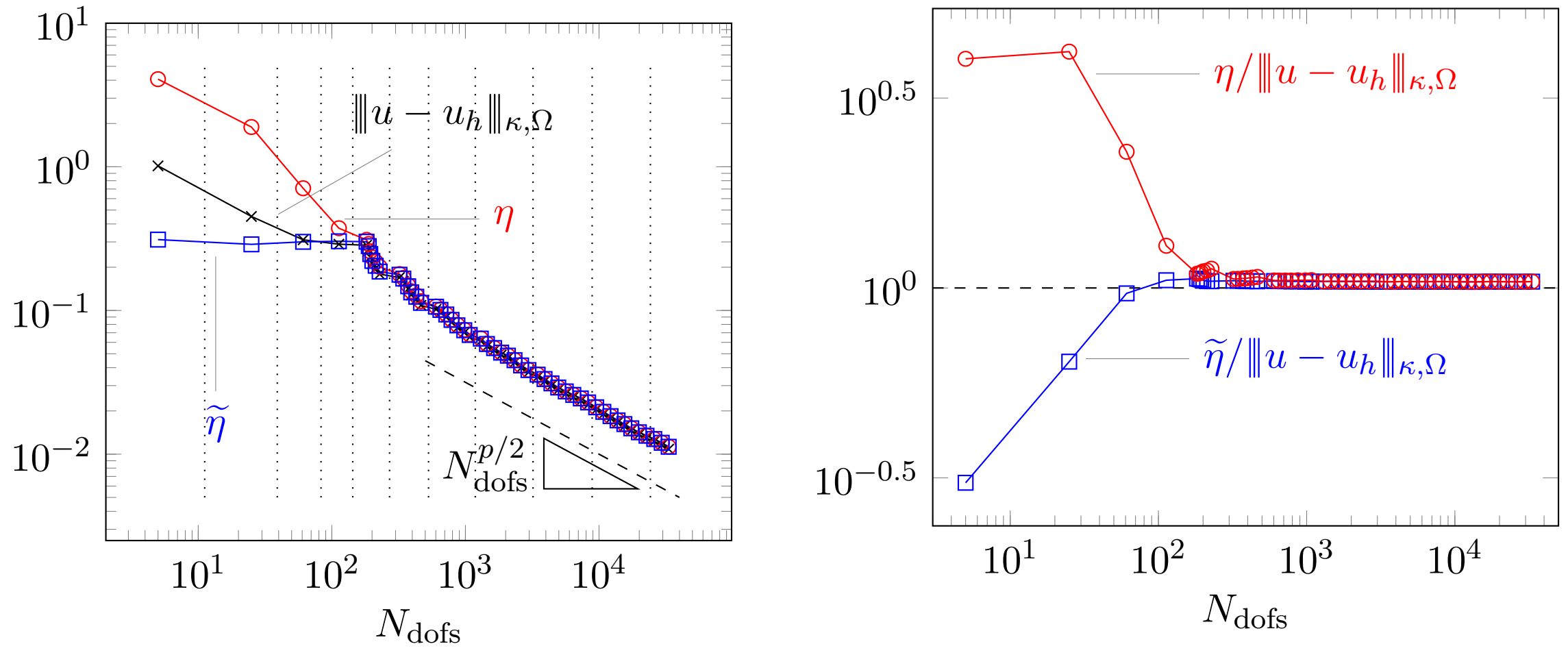

(a) $p = 1$

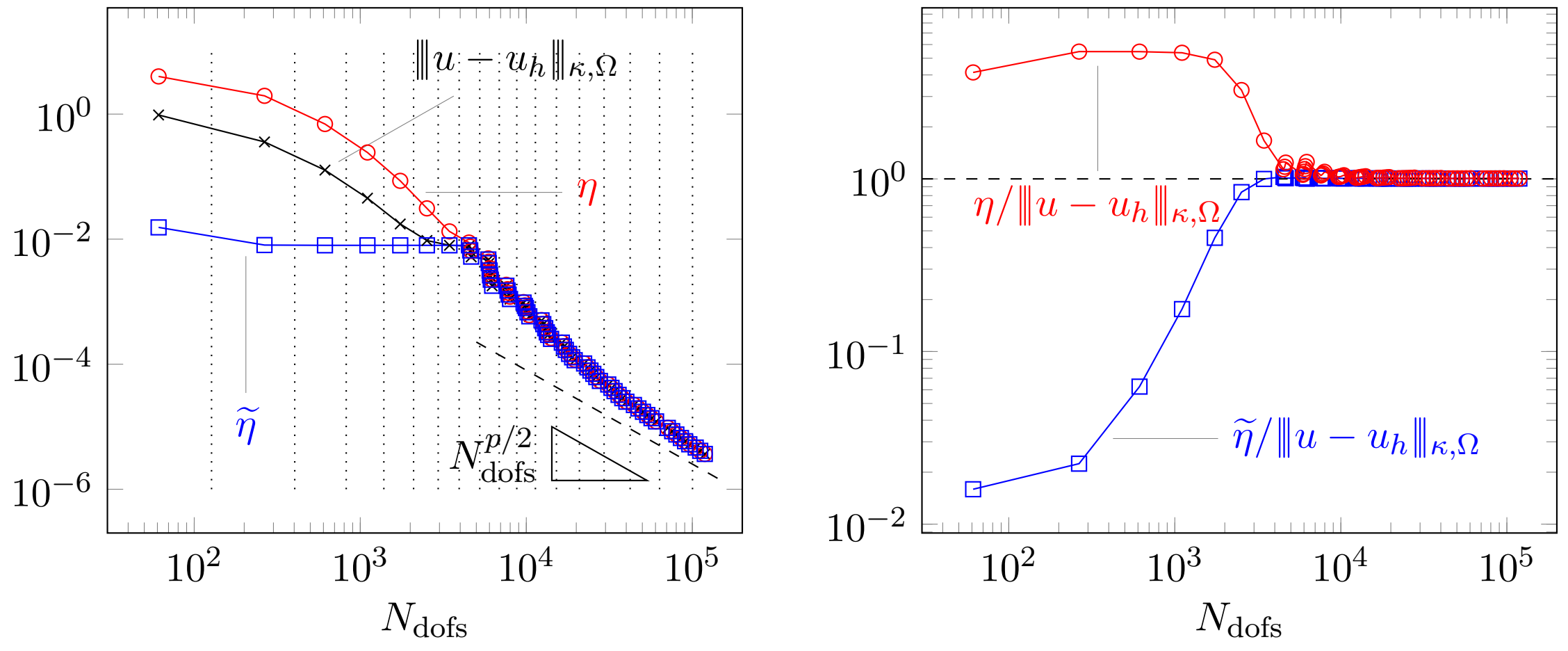

(b) $p = 3$

FIGURE 6.3. Behavior of the error and estimator for the reaction diffusion problem. The dotted vertical lines indicate increases of the truncation parameter $L$ between iterations.

## REFERENCES

1. P.R. Amestoy, I.S. Duff, and J.Y. L'Excellent, *Multifrontal parallel distributed symmetric and unsymmetric solvers*, Comput. Methods Appl. Mech. Engrg. **184** (2000), 501–520.
2. M. Bebendorf, *A note on the Poincaré inequality for convex domains*, Z. Anal. Anwendungen **22** (2003), 751–756.
3. P. Binev, W. Dahmen, and R. De Vore, *Adaptive finite element methods with convergence rates*, Numer. Math. **97** (2004), 219–268.
4. D. Braess, V. Pillwein, and J. Schöberl, *Equilibrated residual error estimates are p-robust*, Comput. Meth. Appl. Mech. Engrg. **198** (2009), 1189–1197.

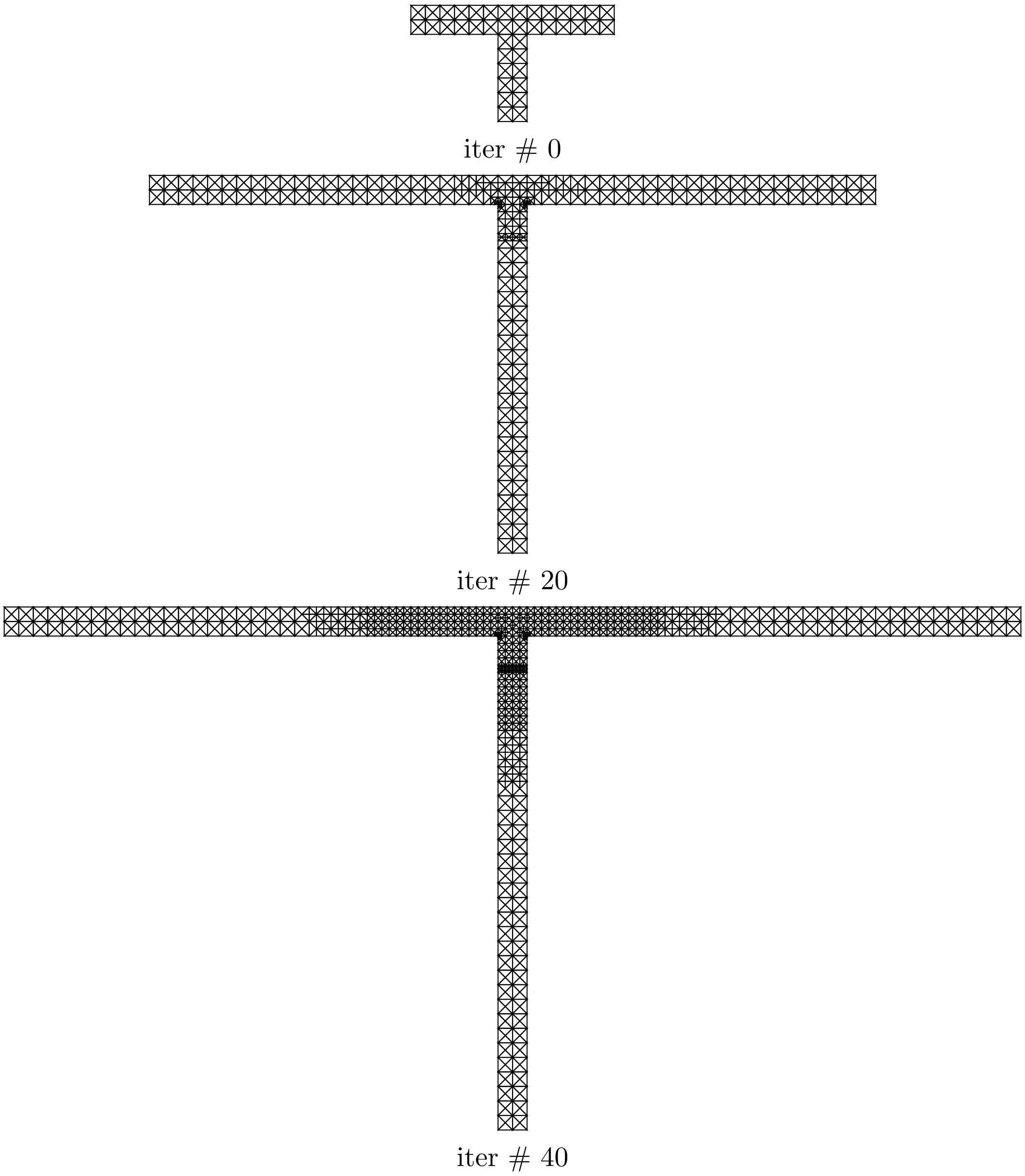

FIGURE 6.4. Sequence of adapted meshes for the Helmholtz problem with $k = 0.7 \cdot 2\pi$

5. T. Chaumont-Frelet, A. Ern, and M. Vohralík, *On the derivation of guaranteed and p-robust a posteriori error estimates for the Helmholtz equation*, Numer. Math. **148** (2021), 525–573.
6. T. Chaumont-Frelet, D. Gallistl, S. Nicaise, and J. Tomezyk, *Wavenumber explicit convergence analysis for finite element discretizations of time-harmonic wave propagation problems with perfectly matched layers*, Comm. Math. Sci. **20** (2022), no. 1, 1–52.
7. T. Chaumont-Frelet and M. Vohralík, *Constrained and unconstrained stable discrete minimizations for p-robust local reconstructions in vertex patches in the de rham complex*, accepted in Found. Comput. Math. (2024).

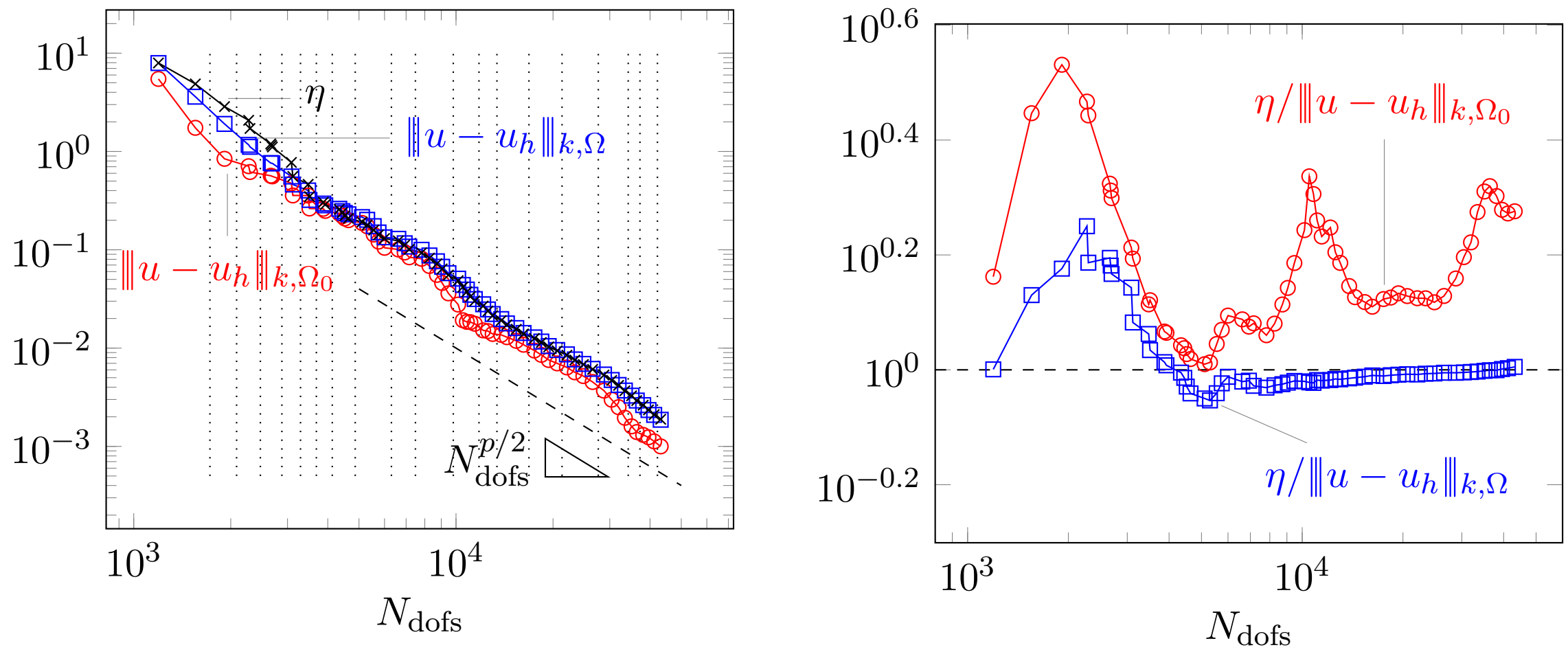

(a) $k = 0.7 \cdot 2\pi$

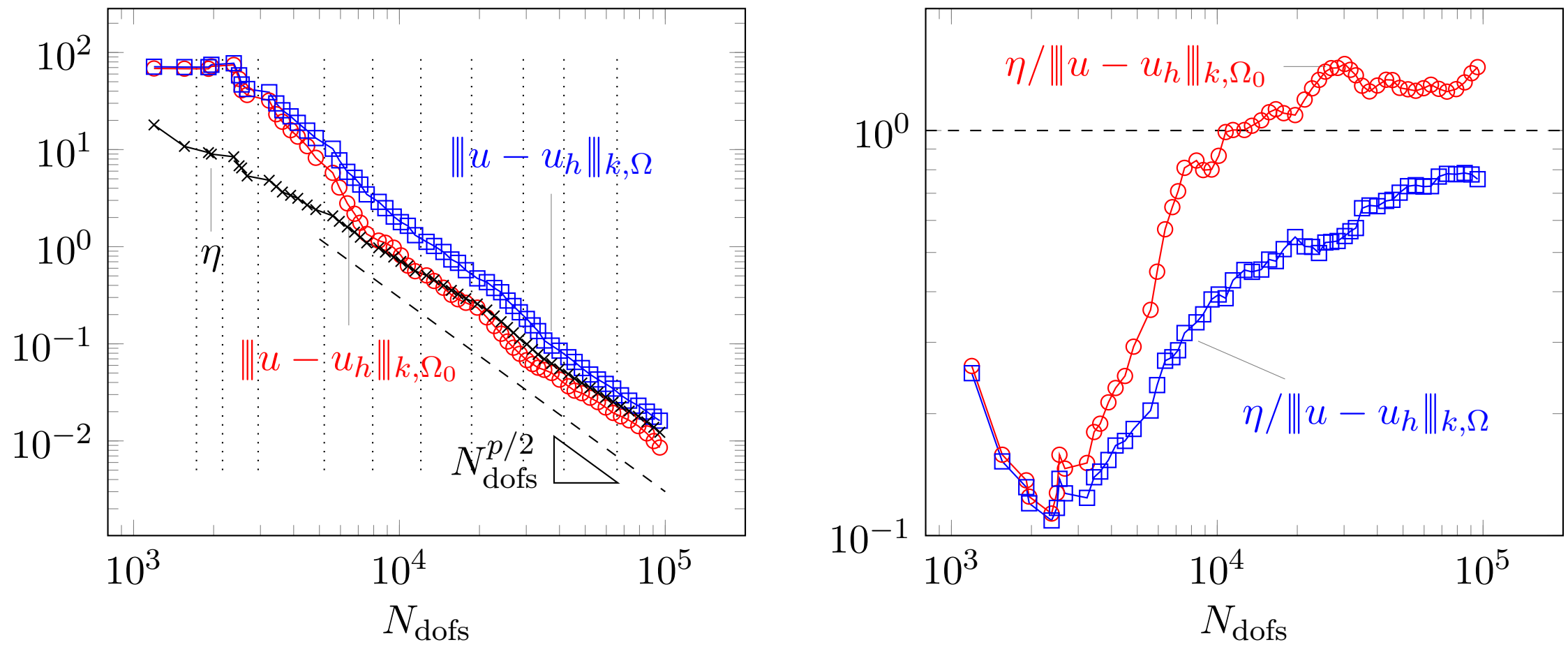

(b) $k = 1.7 \cdot 2\pi$

FIGURE 6.5. Behavior of the error and estimator for the Helmholtz problem. The dotted vertical lines indicate increases of the truncation parameter $L$ between iterations.

8. I. Cheddadi, R. Fučík, M.I. Prieto, and M. Vohralík, *Guaranteed and robust a posteriori error estimates for singularly perturbed reaction-diffusion problems*, ESAIM Math. Model. Numer. Anal. **43** (2009), 867–888.
9. M. Costabel, M. Dauge, and S. Nicaise, *Analytic regularity for linear elliptic systems in polygons and polyhedra*, Math. Models Meth. Appl. Sci. **22** (2012), no. 8, 1250015.
10. P. Destuynder and B. Métivet, *Explicit error bounds in a conforming finite element method*, Math. Comp. **68** (1999), no. 228, 1379–1396.
11. W. Dörfler, *A convergent adaptive algorithm for Poisson's equation*, SIAM J. Numer. Anal. **33** (1996), 1106–1124.
12. W. Dörfler and S. Sauter, *A posteriori error estimation for highly indefinite Helmholtz problems*, Comput. Meth. Appl. Math. **13** (2013), 333–347.
13. A. Ern and J.L. Guermond, *Finite elements I: basic theory and practice*, Springer, 2021.

14. A. Ern and M. Vohralík, *Polynomial-degree-robust a posteriori estimates in a unified setting for conforming, nonconforming, discontinuous Galerkin, and mixed discretizations*, SIAM J. Numer. Anal. **53** (2015), no. 2, 1058–1081.
15. ———, *Stable broken $H^1$ and $\boldsymbol{H}(div)$ polynomial extensions for polynomial-degree-robust potential and flux reconstruction in three space dimensions*, Math. Comp. **89** (2021), 551–594.
16. V. Girault and P.A. Raviart, *Finite element methods for Navier-Stokes equations: theory and algorithms*, Springer-Verlag, 1986.
17. M. Halla, *Analysis fo radial complex scaling methods: scalar resonance problems*, SIAM J. Numer. Anal. **59** (2021), no. 4, 2054–2074.
18. R. Hiptmair and C. Pechstein, *Discrete regular decomposition of tetrahedral discrete 1-forms*, Tech. Report 2017-47, ETH Zürich, 2017.
19. T. Hohage and L. Nannen, *Convergence of infinite element methods for scalar waveguide problems*, BIT Numer. Math. **55** (2015), 215–254.
20. P. Houston, D. Schötzau, and T.P. Wihler, *Energy norm a posteriori error estimation of hp-adaptive discontinuous Galerkin methods for elliptic problems*, Math. Models Meth. Appl. Sci. **17** (2006), no. 1, 33–62.
21. X. Jiang, Z. Sun, L. Sun, and Q. Ma, *An adaptive finite element pml method for helmholtz equations in periodic heterogeneous media*, Comput. App. Math. **43** (2024), no. 4, 242.
22. W. McLean, *Strongly elliptic systems and boundary integral equations*, Cambridge University Press, 2000.
23. J.M. Melenk, *hp-interpolation of nonsmooth functions and an aplication to hp-a posteriori error estimation*, SIAM J. Numer. Anal. **43** (2005), no. 1, 127–155.
24. C. Michler, L. Demkowicz, J. Kurtz, and D. Pardo, *Improving the performance of perfectly matched layers by means of hp-adaptivity*, Numer. Methods Partial Differ. Equ. **23** (2007), no. 4, 832–858.
25. W. Prager and J.L. Synge, *Approximations in elasticity based on the concept of function space*, Quart. Appl. Math. **5** (1947), no. 3, 241–269.
26. C. Schwab, *p- and hp-finite element methods: theory and applications in solid and fluid mechanics*, Clarendon Press, 1998.
27. I. Smears and M. Vohralík, *Simple and robust equilibrated flux a posteriori estimates for singularly perturbed reaction–diffusion problems*, ESAIM: Math. Model. Numer. Anal. **54** (2020), no. 6, 1951–1973.
28. R. Verfürth, *Robust a posteriori error estimators for a singularly perturbed reaction-diffusion equation*, Numer. Math. **78** (1998), 479–493.

## Appendix A. Weighted trace inequality

For completeness, we include a proof of the following standard trace inequality. Other proofs may be found in the literature, e.g., in [8, Eq. (4.4)] and the references therein.

**Lemma A.1** (Weighted trace inequality)**.** *For all $\nu > 0$ and $K \in \mathcal{T}_h$, we have*

$$\rho_K^{-1/2}\|v\|_{\partial K} \leq \max(\beta_K, \sqrt{d+1}(\nu\rho_K)^{-1})\left(\nu^2\|v\|_K^2 + \|\boldsymbol{\nabla} v\|_K^2\right)^{1/2} \tag{A.1}$$

*for all $v \in H^1(K)$.*

*Proof.* We let $\boldsymbol{y}(\boldsymbol{x}) := \boldsymbol{x} - \boldsymbol{x}_K$, where $\boldsymbol{x}_K$ is the center of the larest ball contained in $K$. With $\rho_K$ the radius of the ball, we have $\boldsymbol{y} \cdot \boldsymbol{n}_K \geq \rho_K$ on $\partial K$, leading to

$$\begin{aligned}\rho_K\|v\|_{\partial K}^2 \leq \int_{\partial K}|v|^2\boldsymbol{y}\cdot\boldsymbol{n} = \int_K \boldsymbol{\nabla}\cdot(|v|^2\boldsymbol{y})\\ = \int_K d|v|^2 + \boldsymbol{y}\cdot\boldsymbol{\nabla}|v|^2 = d\|v\|_K^2 + 2\operatorname{Re}\int_K v\boldsymbol{y}\cdot\boldsymbol{\nabla} v\\ \leq d\|v\|_K^2 + 2h_K\|v\|_K\|\boldsymbol{\nabla} v\|_K \leq (d+1)\|v\|_K^2 + h_K^2\|\boldsymbol{\nabla} v\|_K^2.\end{aligned}$$

Division by $\rho_K^2$ gives

$$\rho_K^{-1}\|v\|_{\partial K}^2 \leq \frac{(d+1)}{(\rho_K\nu)^2}\nu^2\|v\|_K^2 + \beta_K^2\|\boldsymbol{\nabla} v\|_K^2,$$

from which (A.1) follows. □

## Appendix B. Radiation condition and perfectly matched layers

In this section, for completeness, we review the radiation condition for our waveguide problem and show that it can be equivalently reformulated with the PML. We then establish that the inf-sup condition in (4.5) holds true under natural assumptions. These results are standard, but we believe it is important to recap them here. Recall that for each $r \in \{1, \dots, R\}$, $\Sigma^r$ is the cross section of the semi-infinite cylinder $C^r_\infty$, and that the $(\lambda^r_j, \psi^r_j)$ are the $L^2(\Sigma^r)$ normalized eigenpairs of the Dirichlet Laplacian on $\Sigma^r$. This notation was rigorously introduced in the beginning of Section 4.

### B.1. Radiation condition.

For $j \geq 0$ and $r \in \{1, \dots, R\}$, the following notation will be useful. If $\lambda^r_j < k$, we let $k^r_j := \sqrt{k^2 - (\lambda^r_j)^2}$ and $k^r_j := i\sqrt{(\lambda^r_j)^2 - k^2}$ otherwise. In all cases $(k^r_j)^2 = k^2 - (\lambda^r_j)^2$.

For $\phi \in H^1_{\rm loc}(\Omega)$ such that $-k^2\phi - \Delta\phi = 0$ in $C_\infty$. We let $\phi^r_\infty = \phi|_{C^r_\infty}$. We then have

$$\phi^r_\infty(\boldsymbol{x}) = \sum_{j \geq 0} \psi^r_j(\boldsymbol{x}_\perp) p^r_j(x_\times)$$

where $p_j$ solves the equation

$$-(k^r_j)^2 p^r_j - (p^r_j)'' = 0.$$

Classically, there then exist two constants $P_\pm$ such that

$$p^r_j(x_\times) = P_- e^{-ik^r_j x_\times} + P_+ e^{ik^r_j x_\times}.$$

We say that $\phi$ satisfies the radiation condition if $P_- = 0$. In other words, we require that for each $1 \leq r \leq R$, there exist constants $\{P^r_j\}_{m \geq 0}$ such that

$$\phi^\infty_r(\boldsymbol{x}) = \sum_{m \geq 0} P^r_j \psi^r_j(\boldsymbol{x}_\perp) e^{ik^r_j x_\times} \tag{B.1}$$

for all $\boldsymbol{x} \in C^r_\infty$.

Given $f : \Omega_0 \to \mathbb{C}$, the original form of the waveguide problem is to find $U : \Omega \to \mathbb{C}$ such that

$$\left\{ \begin{array}{rcll} -k^2 U - \Delta U & = & f & \text{in } \Omega \\ U & = & 0 & \text{on } \partial\Omega \end{array} \right. \tag{B.2}$$

and $U$ satisfies the radiation condition in (B.1).

### B.2. Dirichlet-to-Neumann map.

Let $\Gamma_\infty := \partial C_\infty \setminus \partial\Omega$ be the relatively open surface joining $\Omega_0$ and $C_\infty$. Accordingly, we have $\Gamma_\infty = \cup_{r=1}^R \Gamma^r_\infty$ with $\Gamma^r_\infty = \partial C^r_\infty \setminus \partial\Omega$.

We define the operator $Y^r_\infty : H^{1/2}_{00}(\Gamma^r_\infty) \to H^{-1/2}(\Gamma^r_\infty)$

$$Y^r_\infty \phi = \sum_{m \geq 0} i k^r_m (\phi, \psi^r_m)_{\Gamma^r_\infty} \psi^r_m$$

and correspondingly, $Y_\infty : H^{1/2}_{00}(\Gamma_\infty) \to H^{-1/2}(\Gamma_\infty)$.

B.3. **Truncated Helmholtz problem.** With this notation we can introduce the sesquilinear form

$$b_0(\phi, v) := -k^2(\phi, v)_{\Omega_0} - \langle Y_\infty \phi, v\rangle_{\Gamma_\infty} + (\boldsymbol{\nabla}\phi, \boldsymbol{\nabla} v)_{\Omega_0}$$

for all $\phi, v \in H^1_\Gamma(\Omega_0)$. It is easily seen that we can recast (B.2) into the following weak form: Find $U_0 \in H^1_\Gamma(\Omega_0)$ such that

$$b_0(U_0, v) = (f, v)_{\Omega_0} \tag{B.3}$$

for all $v \in H^1_\Gamma(\Omega_0)$. Then, we have $U_0 = U|_{\Omega_0}$, where $U$ solves (B.2).

In this work, we are interested in wavenumbers $k$ for which the above problem admits a unique solution. Specifically, we will make the assumption that for all $\phi \in H^1_\Gamma(\Omega_0)$

$$b_0(\phi, v) = 0 \quad \forall v \in H^1_\Gamma(\Omega_0) \quad \implies \quad \phi = 0. \tag{B.4}$$

B.4. **Equivalence with the PML problem.** For completeness, we recap the equivalence between the PML problem in (4.2) and the original problem in (B.3).

**Theorem B.1** (Equivalence with the PML problem)**.** *Consider $f \in L^2(\Omega_0)$ and let $u \in H^1_0(\Omega)$ and $U_0 \in H^1_\Gamma(\Omega_0)$ respectively solve (4.2) and (B.3). Then, we have*

$$u|_{\Omega_0} = U_0.$$

*Proof.* For shortness, we write $u_0 = u|_{\Omega_0}$ and we are going to show that $b_0(u_0, v) = (f, v)$ for all $v \in H^1_\Gamma(\Omega_0)$, which lead to the conclusion since the original problem is uniquely solvable by assumption.

We first observe that by construction, $u_0 = 0$ on $\Gamma$, and $-k^2 u_0 - \Delta u_0 = f$ in $\Omega_0$. Therefore, it remains to check that $\boldsymbol{\nabla} u_0 \cdot \boldsymbol{n} = Y_\infty u_0$ on $\Gamma_\infty$, or equivalently, that $\boldsymbol{\nabla} u \cdot \boldsymbol{n} = Y_\infty u$ on $\Gamma_\infty$. We are going to establish this result independently on each $\Gamma^r_\infty$, $1 \leq r \leq R$.

Let us thus fix $r \in \{1, \dots, R\}$. Recalling the definition of $Y^r_\infty$ in Section B.2, introducing the decomposition

$$u(\boldsymbol{x}) = \sum_{j \geq 0} u_j(x_\times)\psi^r_j(\boldsymbol{x}_\perp)$$

we need to show that

$$u'_j(\ell^r) = ik^r_j u_j(\ell^r)$$

for $j \geq 0$. Next, we look at the compatibility condition implied by the fact that $\boldsymbol{\nabla} \cdot (\underline{\boldsymbol{A}}\boldsymbol{\nabla} u) \in L^2(\Omega)$, namely

$$\partial_\times u_j(\ell^r_-) = \alpha^{-1}\partial u_j(\ell^r_+).$$

On the other hand, basic ODE theory shows that

$$u_j(x_\times) = u_j(\ell^r)e^{i\alpha k_r x_\times},$$

from which we conclude that

$$\partial_\times u_j(\ell^r_-) = \alpha^{-1}(i\alpha k^r_j u_j(\ell^r)) = ik^r_j u_j(\ell^r).$$

This conclude the proof. □

B.5. **Well-posedness of the PML problem.** We are now ready to establish the inf-sup stability of $b$.

**Lemma B.2** (Uniqueness inside the PML)**.** *If $\phi \in H^1_0(C_\infty)$ satisfies $b(\phi, v) = 0$ for all $v \in H^1_0(C_\infty)$, then $\phi = 0$.*

*Proof.* We first observe that for $\phi \in H^1_0(C_\infty)$, we have

$$\gamma^{-1} b(\phi,\phi) = -k^2 \|\phi\|^2_{C_\infty} + \|\boldsymbol{\nabla}_\perp \phi\|^2_{C_\infty} + \gamma^{-2} \|\partial_\times \phi\|^2_{C_\infty},$$

since the support of $\phi$ does not intersect $\omega_0$ where $\alpha = 1$. Then, if $\phi$ satisfies the above assumption, we have

$$-k^2 \|\phi\|^2_{C_\infty} + \|\boldsymbol{\nabla}_\perp \phi\|^2_{C_\infty} + \gamma^{-2} \|\partial_\times \phi\|^2_{C_\infty} = 0,$$

and taking the imaginary part reveals that $\partial_\times \phi = 0$. Since $\phi = 0$ on $\Gamma_\infty$, this implies that $\phi = 0$. □

**Theorem B.3** (Inf-sup stability)**.** *There exists a real number $\beta_{\rm st} > 0$ such that*

$$\min_{\substack{\phi \in H^1_0(\Omega) \\ \|\phi\|_{k,\Omega}=1}} \max_{\substack{\phi \in H^1_0(\Omega) \\ \|\phi\|_{k,\Omega}=1}} \operatorname{Re} b(\phi, v) = \frac{1}{\beta_{\rm st}} \tag{B.5}$$

*Proof.* Classically, the inf-sup condition in (B.5) is equivalent to establishing the bound

$$\|u\|_{k,\Omega} \leq \beta_{\rm st} \|\psi\|_{\star,k,\Omega} := \sup_{\substack{v \in H^1_0(\Omega) \\ \|v\|_{k,\Omega}=1}} \langle \psi, v \rangle,$$

whenever $u \in H^1_0(\Omega)$ satisfies

$$b(u, v) = \langle \psi, v \rangle \qquad \forall v \in H^1_0(\Omega)$$

for all $\psi \in H^{-1}(\Omega)$, the dual space of $H^1_0(\Omega)$. Since $b$ satisfies the Gårding inequality in (4.3), we can use Fredholm alternative [22, Theorem 2.27, 2.33 and 2.34], and we simply need check that this bound is satisfied for $\psi = 0$, i.e., uniqueness of solutions.

Let us therefore assume that $u \in H^1_0(\Omega)$ satisfies $b(u, v) = 0$ for all $v \in H^1_0(\Omega)$. As discussed in (B.4) and Theorem B.1 above, for $u_0 := u|_{\Omega_0} \in H^1_\Gamma(\Omega_0)$, we then have

$$b_0(u_0, v) = 0$$

for all $v \in H^1_\Gamma(\Omega_0)$, which, by assumption ensures that $u_0 = 0$. It therefore follows that $u|_{\Omega_0} = 0$, and in particular, $u_\infty := u|_{C_\infty} \in H^1_0(C_\infty)$. From there, we conclude with Lemma B.2 that $u_\infty = 0$. This concludes the proof. □

## Appendix C. Approximation factor

In this section, we give a proof of Theorem 4.5, namely, we establish an upper bound for the approximation $\delta_n$ introduced in (4.17). We first introduce extra notation in Section C.1. We establish some basic Poincaré and trace inequality in Section C.2, and show that the Helmholtz solution decays exponentially in the PML layer in Section C.3. These results are mostly standard, but include their proof here, with explicit constants, for completeness. In Section C.4, we revisit some standard regularity shift result for elliptic PDE. This is required because we are dealing with an unbounded domain with an unbounded boundary. Finally, the actual proof of Theorem 4.5 is presented in Section C.5.

C.1. **Trace norm.** We start by introducing convenient norms to deal with traces. For $U \subset C_\infty$, we will use the following unweighted norm

$$\|v\|^2_{H^1_k(U)} := k^2\|v\|^2_U + \|\boldsymbol{\nabla} v\|^2_U$$

for all $v \in H^1_0(\Omega)$. For $r \in \{1, \dots, R\}$, it is key to observe that

$$\|v\|^2_{H^1_k(C^r_\infty)} = \sum_{j\geq 0} \left\{(k^2 + \lambda_j^2)\|v_j^r\|^2_{(\ell^r,+\infty)} + \|\partial_\times v_j^r\|^2_{(\ell^r,+\infty)}\right\},$$

where the $v_j^r$ are obtained from the spectral decomposition in (4.1). For $q \geq 0$, we also introduce the trace norm

$$\|v\|^2_{H^{1/2}_{00}(\Sigma^r_n)} := \sum_{j\geq 0} \sqrt{k^2 + (\lambda_j^r)^2}|v_j^r(\ell^r + q/k)|^2$$

for $v \in H^1_0(C^r_\infty)$, and if $w \in H^1_0(\Omega)$, we let

$$\|w\|^2_{H^{1/2}_{00}(\Sigma_n)} := \sum_{r=1}^R \|w|_{C^r}\|^2_{H^{1/2}_{00}(\Sigma^r_n)}.$$

C.2. **Basic estimates in cylinders.** In the forthcoming analysis, we will need basic Poincaré and trace inequalities in the PML, i.e., on cylindrical domains. These estimates are classical, but since the proof are short, we include them here for completeness.

**Lemma C.1** (Poincaré inequality)**.** *For all $v \in H^1(\Omega_q)$, we have*

$$k^2\|v\|^2_{C_q} \leq 2qk\|v\|^2_{\Sigma_q} + 4q^2\|\partial_\times v\|^2_{C_q}. \tag{C.1}$$

*Proof.* We establish (C.1) independently in each cylinder $C^r$. Letting $D := \ell^r + q/k$ and using the identity $(|\psi|)' = 2\operatorname{Re}\psi\psi'$, we start with

$$|v(D, \boldsymbol{x}_\perp)|^2 - |v(x_\times, \boldsymbol{x}_\perp)|^2 = 2\operatorname{Re}\int_{s=x_\times}^D v(\boldsymbol{x}_\perp, s)\overline{\partial_\times v(\boldsymbol{x}_\perp, s)}ds$$

which leads to

$$|v(x_\times, \boldsymbol{x}_\perp)|^2 \leq |v(D, \boldsymbol{x}_\perp)|^2 + 2\operatorname{Re}\int_{s=0}^D |v(\boldsymbol{x}_\perp, s)||\partial_\times v(\boldsymbol{x}_\perp, s)|ds,$$

and integration over $\boldsymbol{x}_\perp \in \Sigma$ gives

$$\int_{\boldsymbol{x}_\perp\in\Sigma^r} |v(x_\times, \boldsymbol{x}_\perp)|^2 \leq \|v\|^2_{\Sigma^r_q} + 2\|v\|_{C^r_q}\|\partial_\times v\|_{C^r_q}.$$

We then integrate over $x_\times \in (\ell^r, D)$, leading to

$$\|v\|^2_{C^r_q} \leq D\|v\|^2_{\Sigma^r_q} + 2D\|v\|_{C^r_q}\|\partial_\times v\|_{C^r_q}.$$

The algebraic inequality

$$2D\|v\|_{C^r_q}\|\partial_\times v\|_{C^r_q} \leq \frac{1}{2}\|v\|^2_{C^r_q} + 2D^2\|\partial_\times v\|^2_{C^r_q}$$

then shows that

$$\frac{1}{2}\|v\|^2_{C^r_q} \leq D\|v\|^2_{\Sigma^r_q} + 2D^2\|\partial_\times v\|^2_{C^r_q},$$

and (C.1) follows after multiplying both sides by $2k^2$ and summing over $r \in \{1, \dots, R\}$. $\square$

**Lemma C.2** (Trace inequality). *For $q \geq 0$ and $v \in H^1_0(\Omega)$, we have*

$$k\|v\|^2_{\Sigma_q} \leq \|v\|^2_{H^{1/2}_k(\Sigma_q)} \leq 2\|v\|^2_{H^1_k(C_\infty)} \leq \frac{2|\gamma|^2}{\gamma_{\rm r}}\|v\|^2_{k,\Omega}. \tag{C.2}$$

*Proof.* The first and last inequalities in (C.2) are straightforward, we therefore focus on the middle one.

Let us first establish that for all $\psi \in H^1(\mathbb{R})$ and $x, \tau \geq 0$, we have

$$\tau|\psi(x)|^2 \leq 2\tau^2\|\psi\|^2_{\mathbb{R}_+} + \|\psi'\|^2_{\mathbb{R}_+} \leq 2\left(\tau^2\|\psi\|^2_{\mathbb{R}_+} + \|\psi'\|^2_{\mathbb{R}_+}\right). \tag{C.3}$$

If $\tau = 0$, (C.3) trivially holds true, so that we focus on the case where $\tau > 0$. For $y \geq x$, elementary computations reveal that

$$|\psi(x)|^2 - |\psi(y)|^2 = 2\,\mathrm{Re}\int_x^y \psi(z)\psi'(z)dz \leq 2\|\psi\|_{\mathbb{R}_+}\|\psi'\|_{\mathbb{R}_+}.$$

We then add $|\psi(y)|^2$ on both side, and integrate the resulting inequality in $y$ over $(x, x+1/\tau)$. This gives

$$\frac{1}{\tau}|\psi(x)|^2 \leq \|\psi\|^2_{\mathbb{R}_+} + \frac{2}{\tau}\|\psi\|_{\mathbb{R}_+}\|\psi'\|_{\mathbb{R}_+}.$$

We then obtain (C.3) after multiplying by $\tau^2$ and using Young's inequality on the last term.

Consider now $v \in H^1_0(\Omega)$. For each $r \in \{1, \dots, R\}$, we expend

$$v(x_\times, \boldsymbol{x}_\perp) = \sum_{j\geq 0} v^r_j(x_\times)\psi^r_j(\boldsymbol{x}_\perp)$$

in $C^r$. For each $j \geq 0$, we can then apply (C.3) to $v_j$ with $\tau = \sqrt{k^2 + (\lambda^r_j)^2}$, which gives

$$\|v\|^2_{H^{1/2}_k(\Sigma^r_q)} \leq 2\|v\|^2_{H^1_k(C^r)}$$

after summation over $j$ as per the discussion in Section C.1. The missing estimate in (C.2) then follows by summation over $r$. □

C.3. **Exponential decay in the PML.** We can now establish that the solutions decay inside the PML, and exponentially so. For $j \geq 0$ and $1 \leq r \leq R$, the quantity $\nu^r_j := \mathrm{Re}(\gamma k^r_j)$ will play an important role. The following key properties hold true:

**Lemma C.3** (Properties of $\nu^r_m$). *For all $j \geq 0$ and $r \in \{1, \dots, R\}$, we have*

$$\nu^r_j \geq \gamma_\star|k^r_j| \geq \gamma_\star k_\star > 0 \tag{C.4}$$

*and*

$$\frac{\lambda^r_j}{\nu^r_j} \leq \max(1, \mu)\frac{\sqrt{2}}{\gamma_\star\mu}. \tag{C.5}$$

*Proof.* For shortness, we fix and drop the index $r$ in the notation for this proof. We first establish (C.4). If $k^2_m > 0$, then $\mathrm{Re}(\gamma k_m) = \gamma_{\rm r} k_m = \gamma_{\rm r}|k_m| \geq \gamma_\star k_\star$. On the other hand, when $k^2_m < 0$, we have $k_m = -i|k_m|$, $\gamma k_m = -i\gamma_r|k_m| + \gamma_{\rm i}|k_m|$, and $\mathrm{Re}(\gamma k_m) = \gamma_{\rm i}|k_m| \geq \gamma_\star k_\star$.

To prove (C.5), we first note that

$$\frac{\lambda_m}{\nu_m} \leq \frac{1}{\gamma_\star}\frac{\lambda}{|k_m|}$$

and we then distinguish two cases. First, if $\lambda_m^2 \geq 2k^2$, we have

$$|k_m|^2 = \lambda_m^2 - k^2 \geq \frac{1}{2}\lambda_m^2,$$

leading to

$$\frac{\lambda_m^2}{|k_m|^2} \leq 2.$$

If on the other hand, $\lambda_m^2 \leq 2k^2$, we simply write that

$$\frac{\lambda_m^2}{|k_m|^2} \leq \frac{2k^2}{k_\star^2},$$

and (C.5) follows by taking the square root. □

**Theorem C.4** (Decay in the PML)**.** *Let $0 \leq n \leq m$ and $u \in H^1_{\partial\Sigma}(\Omega_n^{\rm c})$ such that*

$$-k^2\gamma u - \gamma\Delta_\perp u - \gamma^{-1}\partial_\times^2 u = 0 \text{ on } \Omega_n^{\rm c}.$$

*Then, we have*

$$\|u\|_{k,\Omega_m^{\rm c}} \lesssim \frac{e^{-\gamma_\star\mu(m-n)}}{\sqrt{\gamma_\star\mu}}\|u\|_{H_k^{1/2}(\Sigma_n)}, \tag{C.6}$$

*where the hidden constant solely depends on $\nu_{\min}$.*

*Proof.* We are going to establish the two estimates independently in each cylinder $C_\infty^r$, $1 \leq r \leq R$. The estimate then readily follows by summation. Henceforth, we thus consider one fixed cylinder $C_\infty^r$ but we drop the index $r$ to lighten the notation.

Consider the eigenfunction expansions

$$u(\boldsymbol{x}) = \sum_{j\geq 0} u_j(x_\times)\psi_j(\boldsymbol{x}_\perp) \qquad u(\ell^r + n/k, \boldsymbol{x}_\perp) = \sum_{j\geq 0} U_j\psi_j(\boldsymbol{x}_\perp)$$

in $C^r$ and $\Sigma_n^r$. Then, $u_j \in H^1(\ell^r + n/k, +\infty)$ satisfies

$$-(\gamma k_j)^2 u_j - \partial_\times^2 u_j = 0$$

and $u_j(\ell^r + n/k) = U_j$, and basic ODE theory tells us that

$$u_j(x_\times) = U_j e^{-i\gamma k_j(x_\times - n/k - \ell^r)}.$$

and simple computations reveal that

$$|u_j(x_\times)|^2 = |U_j|^2 e^{-2\nu_j(x_\times - n/k - \ell^r)}.$$

Elementary manipulations show that

$$\int_{x_\times=\ell^r+m/k}^{+\infty} e^{-2\nu_j(x_\times - n/k - \ell^r)} dx_\times = \frac{e^{-2\nu_j(m-n)/k}}{2\nu_j},$$

which immediately gives

$$\|u\|^2_{\Omega_m^{\rm c}} = \sum_{j\geq 0}|U_j|^2 \int_{x_\times=\ell^r+m/k}^{+\infty} e^{-2\nu_j(x_\times - n/k)} \leq \frac{e^{-\gamma_\star\mu(m-n)}}{2\gamma_\star k_\star}\|u\|^2_{\Sigma_n},$$

and

$$k\|u\|_{\Omega_m^{\rm c}} \leq \frac{1}{\sqrt{2}}\frac{e^{-\gamma_\star\mu(m-n)}}{\sqrt{\gamma_\star\mu}}\|u\|_{\Sigma_n}. \tag{C.7}$$

The derivatives of $u$ are estimated analogously as follows. On the one hand

$$(\partial_\times u_j)(x_\times) = -\nu_j e^{-\nu_j(x_\times - d)} U_j,$$

and therefore

$$\int_{x_\times=\ell^r+m/k}^{+\infty} |(\partial_\times u_j)(x_\times)|^2 dx_\times = \frac{\nu_m}{2} e^{-2\nu_j(m-n)/k}|U_j|^2 \\ \leq \frac{1}{2}e^{-2\nu_j(m-n)/k}|k_j||U_j|^2 \leq \frac{1}{2}e^{-2\nu_j(m-n)/k}(k+\lambda_j)|U_j|^2,$$

from which we conclude that

$$\|\partial_\times u\|_{\Omega_m^c} \leq \frac{1}{\sqrt{2}} e^{-\gamma_\star \mu(m-n)} \|U\|_{H_k^{1/2}(\Sigma_n)}. \tag{C.8}$$

For the remaining derivatives, we write that

$$\|\boldsymbol{\nabla}_\perp u\|_{\Omega_m^c}^2 = \sum_{j\geq 0} \lambda_j^2 |U_j|^2 \int_{x_\times=\ell^r+m/k}^{+\infty} e^{-2\nu_j(x_\times-d)} dx_\times \\ = \sum_{j\geq 0} \frac{\lambda_j}{\nu_j} \frac{e^{-2\nu_m(m-n)/k}}{2} \lambda_j |U_j|^2 \leq \left( \max_{j\geq 0} \frac{\lambda_j}{\nu_j} \frac{e^{-2\nu_j(m-n)/k}}{2} \right) \|U\|_{H_k^{1/2}(\Sigma_n)}^2, \tag{C.9}$$

and we employ (C.4) and (C.5) to show that

$$\|\boldsymbol{\nabla}_\perp u\|_{\Omega_m^c}^2 \leq \max\left(1, \frac{k}{k_\star}\right) \frac{\sqrt{2}}{\gamma_\star} \frac{e^{-2\gamma_\star(m-n)}}{2} \|U\|_{H_k^{1/2}(\Sigma_n)}^2. \tag{C.10}$$

Now, (C.6) follows by adding (C.7), (C.8) and (C.10). □

**Corollary C.5** (Decay in the PML)**.** *For $v \in H_0^1(\Omega)$ satisfying the assumptions of Theorem C.4, we have*

$$\|v\|_{k,\Omega_m^c} \lesssim \frac{e^{-\gamma_\star \mu(m-n)}}{\sqrt{\gamma_\star \mu}} \|v\|_{k,\Omega}, \tag{C.11}$$

*where the hidden constant depends on $\gamma$ and $\nu_{\min}$.*

*Proof.* We simply combined the trace inequality in (C.2) with (C.6). □

C.4. **Regularity shift.** We next need a regularity shift result. Such estimates are standard on bounded domains [9], but to the best of the author's knowledge, there are no easily citable results in unbounded domains with unbounded boundaries. We therefore establish such a regularity shift estimate in $\Omega$ here.

From here on, for $0 < s < 1$, if $w \in L^2(\Omega)$, we introduce the fractional semi-norm

$$|w|_{H^{1+s}(\Omega)}^2 := \int_\Omega \int_\Omega \frac{|(\boldsymbol{\nabla} w)(\boldsymbol{x}) - (\boldsymbol{\nabla} w)(\boldsymbol{y})|^2}{|\boldsymbol{x}-\boldsymbol{y}|^{d+2s}} d\boldsymbol{y} d\boldsymbol{x},$$

and the corresponding fractional Sobolev space

$$H^{1+s}(\Omega) := \left\{ w \in L^2(\Omega) \quad | \quad |w|_{H^{1+s}(\Omega)} < +\infty \right\}.$$

In the proof below, we employ other fractional Sobolev spaces as well as interpolation inequalities. Our notation are standard, and we refer to [22, Chapter 3 and Appendix B] for a detailed exposition.

**Theorem C.6** (Regularity shift)**.** *Assume that $k \geq k_0 > 0$. Assume that $w \in H_0^1(\Omega)$ satisfies $-\boldsymbol{\nabla} \cdot (\underline{\boldsymbol{A}}\boldsymbol{\nabla} w) \in L^2(\Omega)$. Then, there exists $0 < s < 1$ such that $w \in H^{1+s}(\Omega)$ with*

$$|w|_{H^{1+s}(\Omega)} \lesssim (\|\boldsymbol{\nabla} \cdot (\underline{\boldsymbol{A}}\boldsymbol{\nabla} w)\|_\Omega + k\|w\|_{k,\Omega})^s \|w\|_{k,\Omega}^{1-s} \tag{C.12}$$

*where the hidden constant depends on $k_0$, $s$, $\gamma$ and $\Omega$.*

*Proof.* Throughout this proof, the hidden constants in the $\lesssim$ notation may depend on $k_0$, $s$, $\gamma$ and $\Omega$.

We consider a smooth partition of unity of the real line $(U_m)_{m\geq 0}$ and $U_m(x) = 0$ if $|x-m| \geq 1$. These can be constructed in such a way that

$$|U_m| \leq 1 \qquad |U_m'| \lesssim 1$$

uniformly in $x$ and $m$. We then define a partition of unity for $\Omega$ by setting $\chi_0 = 1$ in $\Omega$ and $\chi_0(\boldsymbol{x}) = U_0(k(x_\times^r - \ell^r))$ in $C_\infty^r$ for $1 \leq r \leq R$, and similarly, for $m \geq 1$, $\chi_m = 0$ in $\Omega$ and $\chi_m(\boldsymbol{x}) = U_m(k(x_\times^r - \ell^r))$ in each $C_\infty^r$. One can readily show that

$$\|\chi_m\|_{L^\infty(\Omega)} \leq 1, \qquad \|\boldsymbol{\nabla}\chi_m\|_{\boldsymbol{L}^\infty(\Omega)} \lesssim k, \qquad \|\Delta\chi_m\|_{L^\infty(\Omega)} \lesssim k^2.$$

We now decompose

$$w = \sum_{m\geq 0} \chi_m w =: \sum_{m\geq 0} w_m.$$

We can readily see that for each $m \geq 0$, we have

$$\boldsymbol{\nabla} \cdot (\underline{\boldsymbol{A}}\boldsymbol{\nabla} w) = \chi_m \boldsymbol{\nabla} \cdot (\underline{\boldsymbol{A}}\boldsymbol{\nabla} w) + 2\underline{\boldsymbol{A}}\boldsymbol{\nabla}\chi_m \cdot \boldsymbol{\nabla} w + \Delta\chi_m w =: f_m$$

with

$$\|f_m\|_\Omega \lesssim \|\boldsymbol{\nabla} \cdot (\underline{\boldsymbol{A}}\boldsymbol{\nabla} w)\|_{V_m} + k\|\boldsymbol{\nabla} w\|_{V_m} + k^2\|w\|_{V_m} \lesssim \|\boldsymbol{\nabla} \cdot (\underline{\boldsymbol{A}}\boldsymbol{\nabla} w)\|_{V_m} + k\|w\|_{k,V_m} \tag{C.13}$$

where we introduced $V_m := \operatorname{supp}\chi_m$. Since each $w_m$ is now compactly supported, we can apply a standard regularity shift in a domain $\widetilde{V}_m \supset V_m$ (see e.g. [9]) to show that

$$|w_m|_{H^{1+s}(\Omega)} \lesssim C_{k,m}\|f_m\|_{H^{s-1}(\Omega)},$$

for some $s > 0$ where $s$ and the constant $C_{k,m}$, in principle, depend on $m, k$ via the shape $V_m$. We will see however that this is not the case, leading to

$$|w_m|_{H^{1+s}(\Omega)} \lesssim \|f_m\|_{H^{s-1}(\Omega)}, \tag{C.14}$$

for a fixed $s$ and hidden constant depending on $\Omega$, $\gamma$ and $k_0$. We consider the case $m = 0$ first. In this case, we simply need to take $\widetilde{V}_0 := \Omega_1$. For $m \geq 1$ we observe each $V_m$ consists of $r$ connected component, each of them being a cylinder of section $\Sigma^r$ and of length $2/k$. We first note that up to a translation (which does not affect the constant in the regularity shift), the connected component of $V_m$ do not depend on $m$, so that the constant in the regularity shift estimate does not depend on $m$. Besides, as for the case $m = 0$, the dependency on $k$ may be lifted by simply considering oversampled domains consisting of cylinders of length $2/k_0$.

We next estimate the right-hand side of (C.14). We start by writing that by interpolation,

$$\|f_m\|_{H^{s-1}(\Omega)} \lesssim \|f_m\|_\Omega^s \|f_m\|_{H^{-1}(\Omega)}^{1-s}. \tag{C.15}$$

We then observe that since $-\boldsymbol{\nabla} \cdot (\underline{\boldsymbol{A}}\boldsymbol{\nabla} w_m) = f_m$, we have

$$|(f_m, v)_\Omega| = |(\underline{\boldsymbol{A}}\boldsymbol{\nabla} w_m, \boldsymbol{\nabla} v)_\Omega| \lesssim \|\boldsymbol{\nabla} w_m\|_\Omega \|\boldsymbol{\nabla} v\|_\Omega,$$

for all $v \in H_0^1(\Omega)$, so that

$$\|f_m\|_{H^{-1}(\Omega)} \lesssim \|\boldsymbol{\nabla} w_m\|_\Omega \lesssim \|w\|_{V_m}, \tag{C.16}$$

where we used the product-rule and the bound on $\boldsymbol{\nabla}\chi_m$ in the last estimate. We now regroup (C.13), (C.15) and (C.16), so that

$$|w_m|_{H^{1+s}(\Omega)} \lesssim \|f_m\|_{H^{s-1}(\Omega)} \lesssim (\|\boldsymbol{\nabla}\cdot(\underline{\boldsymbol{A}}\boldsymbol{\nabla}w)\|_{V_m} + k\|w\|_{V_m})^s \|w\|_{V_m}^{1-s}.$$

We then sum the result leading to

$$\begin{aligned} |w|_{H^{1+s}(\Omega)} &\lesssim \sum_{m\geq 0} (\|\boldsymbol{\nabla}\cdot(\underline{\boldsymbol{A}}\boldsymbol{\nabla}w)\|_{V_m} + k\|w\|_{V_m})^s \|w\|_{V_m}^{1-s}. \\ &\lesssim \left(\sum_{m\geq 0} \|\boldsymbol{\nabla}\cdot(\underline{\boldsymbol{A}}\boldsymbol{\nabla}w)\|_{V_m} + k\|w\|_{V_m}\right)^s \left(\sum_{m\geq 0}\|w\|_{V_m}\right)^{1-s} \\ &\lesssim (\|\boldsymbol{\nabla}\cdot(\underline{\boldsymbol{A}}\boldsymbol{\nabla}w)\|_{\Omega} + k\|w\|_{\Omega})^s \|w\|_{\Omega}^{1-s}, \end{aligned}$$

where we employ Hölder's inequality with $p = 1/s$ and $q = 1/(1-s)$ and then the fact that at most two $V_m$ overlap at any point in $\Omega$. □

C.5. **Approximation factor.** We have now introduced all the necessary tools to estimate the approximation factor $\delta_n$. First, we apply the regularity shift result of Theorem C.6 to solutions of the Helmholtz problem.

**Lemma C.7** (Regularity shift for Helmholtz). *Assume that $k \geq k_0 > 0$, let $\psi \in L^2(\Omega)$ and assume that $\xi \in H^1_0(\Omega)$ solves*

$$b(w,\xi) = (\alpha w, \psi)_\Omega$$

*for all $w \in H^1_0(\Omega)$. Then, $\xi \in H^{1+s}(\Omega)$ with*

$$|\xi|_{H^{1+s}(\Omega)} \lesssim \beta_{\rm st} k^s \|\psi\|_{\alpha_{\rm r},\Omega}. \tag{C.17}$$

*In addition, if $n \geq 0$, $m \geq n+1$ and $\psi \in L^2(\Omega_n)$, there exists a smooth cutoff function $\chi$ such that $\chi = 1$ on $\Omega^{\rm c}_{m+1}$, $\chi = 0$ on $\Omega_m$,*

$$\|\chi\|_{L^\infty(\Omega)} \leq 1 \qquad \|\boldsymbol{\nabla}\chi\|_{\boldsymbol{L}^\infty(\Omega)} \lesssim k \qquad \|\Delta\chi\|_{L^\infty(\Omega)} \lesssim k^2$$

*and*

$$\|\chi\xi\|_{k,\Omega} + k^{-s}|\chi\xi|_{H^{1+s}(\Omega)} \lesssim \beta_{\rm st}\frac{e^{-\gamma_\star\mu(m-n)}}{\sqrt{\gamma_\star\mu}}\|\psi\|_{\alpha_{\rm r},\Omega} \tag{C.18}$$

*holds true. The hidden constants in (C.17) and (C.18) depend on $k_0$, $s$, $\gamma$, and $\Omega$.*

*Proof.* Throughout this proof, the hidden constants in the $\lesssim$ notation are allowed to depend on $k_0$, $s$, $\gamma$ and $\Omega$.

The estimate in (C.17) is a direct consequence of (C.12), Indeed, we have on the one hand

$$\operatorname{Re} b(\xi, v) = \operatorname{Re} k(\alpha\psi, v) \lesssim \|\psi\|_{\alpha_{\rm r},\Omega}\|v\|_{k,\Omega}$$

for all $v \in H^1_0(\Omega)$, and it then follows from (B.5) that

$$\|\xi\|_{k,\Omega} \lesssim \beta_{\rm st}\|\psi\|_{\alpha_{\rm r},\Omega} \lesssim (1+\beta_{\rm st})\|\psi\|_{\alpha_{\rm r},\Omega}.$$

On the other hand, we have

$$\|\boldsymbol{\nabla}\cdot(\underline{\boldsymbol{A}}\boldsymbol{\nabla}\xi)\|_\Omega = \|k\alpha\psi + k^2\alpha\xi\|_\Omega \lesssim k\|\psi\|_{\alpha_{\rm r},\Omega} + k\|\xi\|_{k,\Omega} \lesssim (1+\beta_{\rm st})k\|\psi\|_{\alpha_{\rm r},\Omega},$$

which concludes the proof of (C.17). Indeed, $\beta_{\rm st} \gtrsim 1$, see [6, Lemma 3.2].

We now turn to (C.18). Since $\operatorname{supp}\chi \cap \operatorname{supp}\psi = \emptyset$, we have

$$b(\chi\xi, v) = (2\underline{\boldsymbol{A}}\boldsymbol{\nabla}\chi\cdot\boldsymbol{\nabla}\xi + \Delta\chi\xi, v)_\Omega,$$

and therefore

$$\|\boldsymbol{\nabla} \cdot (\underline{\boldsymbol{A}}(\chi\xi))\|_\Omega \leq \|k^2\alpha\chi\xi + 2\underline{\boldsymbol{A}}\boldsymbol{\nabla}\chi \cdot \boldsymbol{\nabla}\xi + \Delta\chi\xi\|_\Omega \lesssim k\|\xi\|_{\Omega_m^c}.$$

Similarly, we have

$$\|\chi\xi\|_\Omega \lesssim \|\xi\|_{\Omega_m^c},$$

and we conclude from (C.12) that

$$|\chi\xi|_{H^{1+s}(\Omega)} \lesssim k^s \|\xi\|_{k,\Omega_m^c}.$$

Now (C.18) follows from (C.11) since

$$\|\xi\|_{k,\Omega_m^c} \lesssim \frac{e^{-\gamma_\star\mu(m-n)}}{\sqrt{\gamma_\star\mu}} \|\xi\|_{k,\Omega} \leq \beta_{\rm st} \frac{e^{-\gamma_\star\mu(m-n)}}{\sqrt{\gamma_\star\mu}} \|\psi\|_{\alpha_{\rm r},\Omega}.$$

□

We are now ready to bound the approximation factor.

*Proof of Theorem 4.5.* Here, we will need the quasi-interpolation $J_h : H_0^1(\Omega^h) \to V_h$ introduced in [13, Theorem 19.14]. The crucial property we need is that

$$\|v - J_h v\|_{\Omega^h}^2 \preceq_{\mathcal{T}_h} \sum_{K\in\mathcal{T}_h} h_K^2 \|\boldsymbol{\nabla} v\|_K^2, \qquad \|\boldsymbol{\nabla}(v - J_h v)\|_{\Omega^h}^2 \preceq_{\mathcal{T}_h} \sum_{K\in\mathcal{T}_h} h_K^{2s} |v|_{H^{1+s}(K)}^2,$$

which implies in particular that

$$\begin{aligned}\|v - J_h v\|_{\Omega^h}^2 \preceq_{\mathcal{T}_h} &\sum_{K\in\mathcal{T}_h} \left( h_K^2 \|\boldsymbol{\nabla} v\|_K^2 + h_K^{2s} |v|_{H^{1+s}(K)}^2 \right) \\ &\leq (kh_m)^2 \|\boldsymbol{\nabla} v\|_\Omega^2 + (kh_m)^{2s} k^{-2s} |v|_{H^{1+s}(\Omega)}^2\end{aligned}$$

for all $v \in H_0^1(\Omega^h) \cap H^{1+s}(\Omega)$.

In the remainder of this proof, the hidden constants in the $\lesssim$ notation are allowed to depend on $k_0$, $s$, $\gamma$ and $\Omega$, and $\beta_{\mathcal{T}_h}$.

We now consider a right-hand side $\psi \in L^2(\Omega_n)$ and consider the associated (adjoint) solution $\xi \in H_0^1(\Omega)$. We introduce $\widetilde{\xi} := \xi - \chi\xi$, where $\chi$ is the cutoff function from Lemma C.7. Then, $\widetilde{\xi} \in H_0^1(\Omega) \cap H^s(\Omega)$ with $\operatorname{supp} \widetilde{\xi} \subset H_0^1(\Omega_{m+1})$ and

$$\|\xi - \widetilde{\xi}\|_{k,\Omega} \lesssim \beta_{\rm st} \frac{e^{-\gamma_\star\mu(m-n)}}{\sqrt{\gamma_\star\mu}} \|\phi\|_{\alpha_{\rm r},\Omega}, \qquad k^{-s} |\widetilde{\xi}|_{H^s(\Omega)} \lesssim \beta_{\rm st} \|\phi\|_{\alpha_{\rm r},\Omega}.$$

We can now write

$$\begin{aligned}\min_{\xi_h\in V_h} \|\xi - \xi_h\|_{k,\Omega} \leq \|\xi - J_h\widetilde{\xi}\|_{k,\Omega} &\leq \|\xi - \widetilde{\xi}\|_{k,\Omega} + \|\widetilde{\xi} - J_h\widetilde{\xi}\|_{k,\Omega} \\ &\lesssim \|\xi - \widetilde{\xi}\|_{k,\Omega} + (kh_m)^s k^{-s} |\widetilde{\xi}|_{H^s(\Omega)} + kh_m \|\boldsymbol{\nabla}\widetilde{\xi}\|_\Omega.\end{aligned}$$

It follows that

$$\min_{xi_h\in V_h} \|\xi - \xi_h\|_{k,\Omega} \lesssim \beta_{\rm st} \left( \frac{e^{-\gamma_\star\mu(m-n)}}{\sqrt{\gamma_\star\mu}} + (kh_m)^s + kh_m \right) \|\psi\|_{\alpha_{\rm r},\Omega},$$

which completes the proof. □